\RequirePackage{luatex85}
\documentclass[12pt]{amsart}
\usepackage{natbib}
\usepackage{bibunits}
\usepackage[english]{babel}
\usepackage[T1]{fontenc}
\usepackage[utf8]{inputenc}
\usepackage{amsmath,amssymb,amsthm}
\usepackage{textcomp}
\usepackage{stmaryrd}
\usepackage[sc]{mathpazo}
\usepackage[left=2cm,right=2cm,top=3cm,bottom=3cm]{geometry}
\usepackage{eulervm}
\usepackage[cal=dutchcal,scr=kp]{mathalpha}
\makeatletter
\def\mathalfa@frakscaled{s*[1.0]}
  \DeclareFontFamily{U}{esstixfrak}{\skewchar \font =45}
  \DeclareFontShape{U}{esstixfrak}{m}{n}{
    <-> \mathalfa@frakscaled esstixfrak}{}
  \DeclareMathAlphabet{\mathfrak}{U}{esstixfrak}{m}{n}
\makeatother

\usepackage{mathtools}
\usepackage{comment}

\usepackage{tikz}
\usetikzlibrary{knots,hobby,decorations.pathreplacing,shapes.geometric,calc,decorations.text,shapes.misc,calc,math,arrows,decorations.markings,patterns, matrix, intersections,cd}

\usepackage[colorlinks=true,linkcolor=blue]{hyperref}

\definecolor{couleur3}{rgb}{0.0, 0.5, 1.0} 
\definecolor{couleur4}{rgb}{0.55, 0.71, 0.0} 
\definecolor{couleur1}{rgb}{0.87, 0.45, 1.0} 
\definecolor{couleur2}{rgb}{0.81, 0.06, 0.13}
\definecolor{couleur5}{rgb}{1, .64, 0.28}

\newcommand{\sitename}[1]{\mathfrak{#1}}
\newcommand{\nom}[1]{\textsc{#1}}

\newcommand{\C}{\mathbb{C}}
\newcommand{\D}{\mathbb{D}}

\newcommand{\N}{\mathbb{N}}

\renewcommand{\P}{\mathbb{P}}
\newcommand{\Q}{\mathbb{Q}}
\newcommand{\R}{\mathbb{R}}
\renewcommand{\S}{\mathbb{S}}

\newcommand{\Z}{\mathbb{Z}}

\newcommand{\fcal}{\mathcal{f}}

\newcommand{\Ical}{\mathcal{I}}

\newcommand{\Mcal}{\mathcal{M}}

\newcommand{\Ocal}{\mathcal{O}}

\newcommand{\Bscr}{\mathscr{B}}

\newcommand{\Dscr}{\mathscr{D}}

\newcommand{\Sscr}{\mathscr{S}}
\newcommand{\Tscr}{\mathscr{T}}
\newcommand{\Uscr}{\mathscr{U}}
\newcommand{\Vscr}{\mathscr{V}}

\newcommand{\Xscr}{\mathscr{X}}

\newcommand{\Sfrak}{\mathfrak{S}}

\renewcommand{\hbar}{\overline{h}}

\newcommand{\lcal}{\mathcal{l}}

\DeclareMathOperator{\Hom}{\mathrm{Hom}}

\DeclareMathOperator{\Vect}{\mathrm{Vect}}

\DeclareMathOperator{\GL}{\mathrm{GL}}

\DeclareMathOperator{\Cone}{\mathrm{Cone}}

\DeclareMathOperator{\Card}{\mathrm{Card}}
\DeclareMathOperator{\Conv}{\mathrm{Conv}}

\usepackage{cancel}

\begin{document}
\title{Birational morphisms in quantum toric geometry}
\author{Antoine \textsc{Boivin}}
\date\today

\email{\href{mailto:antoine.boivin@univ-angers.fr}{antoine.boivin@univ-angers.fr}}


\begin{abstract}
In this paper, we investigate birational toric morphisms between quantum toric stacks---namely, toric (analytic) stacks associated with fans whose cones may be irrational---focusing on two primary classes of examples: weighted blow-ups with arbitrary weights, and morphisms induced by cobordisms.
\end{abstract}
\maketitle

\newtheorem{Thm}{Theorem}[subsection]

\newtheorem{Prop}[Thm]{Proposition}

\newtheorem{Lemme}[Thm]{Lemma}

\newtheorem{Cor}[Thm]{Corollary}

\newtheorem{Const}[Thm]{Construction}

\newtheorem{DefProp}[Thm]{Definition-Proposition}

\theoremstyle{definition}

\newtheorem{Ex}[Thm]{Example}

\newtheorem{Def}[Thm]{Definition}

\newtheorem{Not}[Thm]{Notation}

\theoremstyle{remark}
\newtheorem{Rem}[Thm]{Remark}
\newtheorem{Avert}[Thm]{Warning}

\setcounter{secnumdepth}{4}
\setcounter{tocdepth}{3}
\tableofcontents

\section*{Introduction}

The combinatorial description of toric varieties --- by fans of rational strongly convex cones of $\R^d$ for toric varieties --- makes toric geometry fertile ground for testing conjectures on general varieties. This is notably true in birational geometry: for instance, the proof of the weak factorization of birational morphisms proved in \cite{abramovich2002torification} uses a reduction to the toric case proved in \cite{wlodarczyk1999birational}. Another example is the Minimal Model Program (MMP): the general higher-dimensional MMP is still conjectural, whereas the toric MMP is fully established (see \cite{Reid1983} for the complete toric varieties case and \cite{fujino2004introduction} and \cite{fujino2006equivariant} for the general case). We can extend this description in order to obtain (algebraic) toric stacks (see \cite{Borisov_2004,fantechi2009smooth,Toricstacks}) which have also been studied in birational geometry (see for instance \cite{Levchenko2022BirationalIO}, \cite{schmitt2024birational}) and, in particular, via weighted blow-ups: for instance, the \nom{Oda} conjecture \cite{Oda1978LecturesOT} of factorization of proper birational maps by blow-ups is true for 3-dimensional toric stacks and weighted blow-ups by \cite{ewald1986stellare} (see \cite{Abramovich2023}).

This description also leads to a link with LVM manifolds. They are a family of non-Kähler compact complex manifolds, introduced in \cite{de1997new}, then generalized in \cite{meers1998}, built from points of $\C^m$ ($m=1$ for the first paper and arbitrary $m$ for the second paper) fulfilling combinatorial conditions. In \cite{LV}, the authors prove that, under some rationality condition, these manifolds are fiber bundles over projective orbifold toric varieties. Moreover, we can describe holomorphic surgeries between two LVM manifolds with different combinatorics of points (see \cite{bosio2006}). These surgeries allow us to describe bimeromorphic maps between these two manifolds which descend to a birational map between the associated toric varieties (see also \cite{pasquier2012approachminimalmodelprogram} for horospherical varieties).

Both examples are generalized by quantum toric stacks. These toric objects, introduced in \cite{katzarkov2020quantum}, are described by fans of strongly convex cones generated by elements of a (not necessarily discrete) subgroup of $\R^d$. Moreover, they are the base space of a $\C^{m}$-fiber bundle whose total space is a LVMB manifold (a generalization of LVM manifolds introduced in \cite{bosio2001}). One of the main features of these toric stacks is the existence of continuous families of quantum toric stacks induced by continuous families of fans. In particular, we get moduli spaces of such toric stacks (studied in \cite{boivin2023moduli} and \cite{boivin2025gluingmodulispacesquantum}).

The goal of this paper is to study generalizations of the blow-up of toric varieties and the surgery of LVM manifolds within the framework of quantum toric stacks. For this purpose, we introduce a suitable definition of birational toric morphism of quantum toric stacks which encompasses the two previous examples.

Firstly, we study the simplest example of blow-up in the quantum case: the (weighted) blow-up of the quantum affine plane. Two main cases appear: the rational weight case, where there exists a morphism between the blow-up and the plane (up to a covering in the non-integer case), and the irrational case, where the blow-up morphism is well-defined on an open dense substack.

We generalize this toy case in two ways. In section \ref{section : rationalwbu}, we will consider the case where the "natural blow-up" morphism is well-defined (i.e., when we add an $\N$-linear combination of the generators of the cones), and we will study the associated exceptional divisor (thanks to computations of the fiber of toric morphisms in subsection \ref{subsection : fibers}). 
In section \ref{general_blowup}, we consider the second point of view, i.e., the birational one. We develop the notion of a birational fan morphism, which corresponds to a torus-equivariant birational morphism between the associated quantum toric stacks. This permits us to understand the irrational blow-up and, in particular, the exceptional divisor of such a blow-up.

In conclusion of this first part of the paper, we prove in section \ref{9-Decomp_Def : morph_birat} that any pair of quantum fans with the same support can be linked by a zig-zag of blow-ups, generalizing a statement of \cite{10.2307/2155318} to the quantum case.

The second part is devoted to the study of cobordisms. In subsection \ref{10-Polytopes_quadriques}, we review the rational case as described in \cite{bosio2006} and the link between the cobordism of the polytopes and the (holomorphic) surgery of associated LVM manifolds. We generalize in section \ref{fan_cobord} the notion of cobordism between quantum fans (which corresponds to the notion between polytopes thanks to the construction of the normal fan). This construction induces a birational morphism between the induced quantum toric stacks (a blow-up of any weight is an example of such a birational morphism). This notion is stable under deformation, and moreover, a cobordism induces a family of quantum toric stacks (by considering slices of the cobordism), which is not the case in the rational framework (since generically a slice is not a rational polytope).

\textbf{Acknowledgement :} The author was supported by the Région Pays de la Loire from the France 2030 program, Centre Henri Lebesgue ANR-11-LABX-0020-01.

\section{Quantum toric stacks}

In this section, we recall the needed definitions and theorems on quantum toric stacks (see \cite{katzarkov2020quantum} and \cite{boivin2020nonsimplicial} for the details of the constructions) and, in particular, their combinatorial description and how to describe the (equivariant) morphisms between them. From this description, we explicitly compute the fiber over an orbit of the torus by a toric morphism between (affine) quantum toric stacks.

\subsection{Recollection on quantum toric stacks}

In the same way as in classical toric geometry, we begin the definition of adapted combinatorial objects in order to describe the considered geometric objects.

\begin{Def}
Let $\Gamma$ be a finitely generated subgroup of $\R^d$ such that $\Vect_\R(\Gamma)=\R^d$. A calibration of $\Gamma$ is given by:
\begin{itemize}
\item A group epimorphism $h \colon \Z^n \to \Gamma$ ;
\item A subset $\Ical \subset \{1,\ldots,n\}$ such that $\Vect_\C(h(e_j), j \notin \Ical)=\C^d$ (this is the set of "virtual generators").

\end{itemize}
This is a standard calibration\footnote{up to an isomorphism, it is always the case} if $\Z^d \subset \Gamma$, $h(e_i)=e_i$ for $i=1,\ldots,d$ and $\Ical$ is of the form $\{n-|\Ical|+1,\ldots,n\}$.
\end{Def}

We use this description of $\Gamma$ since, in the irrational setting, the subgroups of $\R^d$ do not have a fixed rank: the rank of the group $\Z^2+\alpha \Z$ is $2$ if $\alpha \in \Q$ and is $3$ otherwise. This is used in the study of moduli spaces of quantum toric stacks (see \cite[section 11]{katzarkov2020quantum} and \cite{boivin2023moduli}).

\begin{Def} \label{calib_q_fan_def}
A quantum fan $(\Delta,h : \Z^n \to \Gamma \subset \R^d,\Ical)$ in $\Gamma$ is the data of 
\begin{itemize}
\item a collection $\Delta$ of strongly convex polyhedral cones generated by elements of $\Gamma$ such that every intersection of cones of $\Delta$ is a cone of $\Delta$, every face of a cone of $\Delta$ is a cone ;
\item a standard calibration $h$ with $\Ical$ its set of virtual generators;
\item A set of generators $A$ i.e. a subset of $\{1,\ldots,n\} \setminus \Ical$ such that the 1-cone generated by the $h(e_i)$ for $i \in A$ are exactly the 1-cones of $\Delta$.
\end{itemize}  
\begin{Rem}
This is an "irrational" generalization of the extended stacky fan of \cite{jiang2008}.
\end{Rem}

The fan is said simplicial if every cone of $\Delta$ is simplicial (i.e. which can be send on a cone $\Cone(e_1,\ldots,e_k)$ by a linear automorphism of $\R^d$). \\
We note $\Delta(1)$ the set of cones of dimension 1 of $\Delta$ and $\Delta_{max}$ the set of maximal cones of $\Delta$. \\ 
With linear morphisms which preserve inclusion of cones and calibration, the quantum fans form a category denoted $\mathbf{QF}$. 
\end{Def}

\begin{Def}
    A subfan of a quantum fan $(\Delta,h,\Ical)$ is a 
    subobject of $(\Delta,h,\Ical)$ in $\mathbf{QF}$ i.e. a
    quantum fan $(\Delta',h,\Ical')$ such that $\Delta'\subset \Delta$ and $\Ical \subset \Ical'$.
\end{Def}

In \cite{katzarkov2020quantum} and \cite{boivin2020nonsimplicial}, the authors give the construction of a quantum toric stack associated to a quantum fan (simplicial fan in the former and arbitrary fan in the latter).

They are stacks over the equivariant analytic site \footnote{i.e. the category of analytic spaces with an action of an abelian Lie group with a dense orbit, equivariant holomorphic morphisms and coverings by equivariant Euclidean open subsets} $\sitename{A}$ with a "quantum torus" which is dense in it.
The quantum torus $\Tscr_{h,\Ical}$ is the quotient  of $(\C^*)^d$ by the action of $\Z^ {n-d}$ through the map 
\[
x \in \Z^{n-d} \mapsto E(h(0_{\Z^{d}}\oplus x)) \in (\C^*)^d
\]
where $E$ is the map \begin{equation} \label{def_exp}
    (z_1,\ldots,z_n) \in \C^d \mapsto (\exp(2i\pi z_1),\ldots,\exp(2i\pi z_n)) \in (\C^*)^d
\end{equation}
and $(h,\Ical)$ is a standard calibration. This stack is isomorphic to the quotient stack $[\C^d/\Z^n]$ --- the action is given by the morphism $h$ --- thanks to the exponential morphism. Note that if $\mathrm{im}(h) \subset \Z^d$ then $(\C^*)^d$ is a $\Z^{n-d}$-gerbe over $\Tscr_{h,\Ical}$. Such stack is a Picard stack (cf. \cite[exposé XVIII définition 1.4.5]{SGA4}) by descent of the group law of $(\C^*)^d$.

A torus morphism $\Tscr_{h,\Ical} \to \Tscr_{h',\Ical'}$ is a Picard stack morphism $\mathcal{l}\colon\Tscr_{h,\Ical} \to \Tscr_{h',\Ical'}$ which respects combinatorial conditions (\cite[Definition 4.9]{katzarkov2020quantum}): more concretely it is combinatorially described by the data of two linear morphisms $(L \colon \R^d \to  \R^{d'},H : \R^n \to \R^{n'})$ which respects virtual generators and such that the following diagram
\begin{align}\label{torus_morphism}
\begin{tikzcd}[ampersand replacement=\&]
	{\R^n} \& {\R^{n'}} \\
	{\R^d} \& {\R^{d'}}
	\arrow["H", from=1-1, to=1-2]
	\arrow["h"', from=1-1, to=2-1]
	\arrow["{h'}", from=1-2, to=2-2]
	\arrow["L"', from=2-1, to=2-2]
\end{tikzcd}
\end{align}
commutes (see \cite[Proposition 3.3.5]{Boivin}). Indeed, thanks to these conditions, the morphism $L_\C \coloneqq L \otimes id_\C \colon \C^d \to \C^{d'}$ descends to a morphism $\Tscr_{h,\Ical} \to \Tscr_{h',\Ical'}$:
\begin{align}\label{induced_torus_morph}
\begin{tikzcd}[ampersand replacement=\&]
	{\C^d} \& {\C^{d'}} \\
	{\Tscr_{h,\Ical}} \& {\Tscr_{h',\Ical'}}
	\arrow["{L_\C}", from=1-1, to=1-2]
	\arrow[two heads, from=1-1, to=2-1]
	\arrow[two heads, from=1-2, to=2-2]
	\arrow["\lcal",from=2-1, to=2-2]
\end{tikzcd}
\end{align}
This definition permits to define torus morphism "with irrational exponents". 
\begin{Ex}
     \[\begin{tikzcd}[scale=2,ampersand replacement=\&]
	{\Z^2} \&\& {\Z^2} \\
	{\Z+\sqrt{2}\Z} \&\& {\Z+\sqrt{2}\Z} \\
	\C \&\& \C \\
	{\Tscr_{h}=[\C^*/\Z]} \&\& {\Tscr_{h}=[\C^*/\Z]}
	\arrow["{(x,y) \mapsto (y,2x)}", from=1-1, to=1-3]
	\arrow["{z \mapsto z\sqrt{2}}", from=2-1, to=2-3]
	\arrow["{\left\langle \--,(1,\sqrt{2}) \right\rangle}"', from=1-1, to=2-1]
	\arrow["{\left\langle \--,(1,\sqrt{2}) \right\rangle}", from=1-3, to=2-3]
	\arrow[hook,from=2-1, to=3-1]
	\arrow[hook,from=2-3, to=3-3]
	\arrow["{[E]}",from=3-1, to=4-1]
	\arrow["{\text{«}z \mapsto z^{\sqrt{2}}\text{»}}", from=4-1, to=4-3]
	\arrow["{z \mapsto z\sqrt{2}}", from=3-1, to=3-3]
	\arrow["{[E]}",from=3-3, to=4-3]
\end{tikzcd}\]
where $\Z$ acts on $\C^*$ through
\[
n \cdot t=E(n\sqrt{2}) t.
\]
\end{Ex}

Let $(\Delta,h,\Ical)$ a quantum fan. The construction of the quantum toric stack associated to it is done as follows:

Every simplicial cone $\sigma=\Cone(h(e_{i_1}),\ldots,h(e_{i_{\dim(\sigma)}}))$ (isomorphic to $\Cone(e_1,\ldots,e_{\dim(\sigma)})$ by a linear morphism $L_\sigma$ and note $H$ a linear morphism described by a permutation $\chi$ of $\{1,\ldots,n\}$ such that $\chi(i_k)=k$ for $1 \leq k \leq \dim(\sigma)$) of a simplicial quantum fan $(\Delta,h,\Ical)$ describe an affine quantum toric stack
\[
\Uscr_\sigma \coloneqq [\C^{\dim(\sigma)}\times (\C^*)^{d-\dim(\sigma)}/\Z^{n-d}]
\]
where the action of $\Z^{n-d}$ on $\C^{\dim(\sigma)}\times (\C^*)^{d-\dim(\sigma)}$ is given by the following morphism 
\[
x \in \Z^{n-d} \mapsto  EL_\sigma h H^{-1}(0_{\R^{d}} \oplus x) \in (\C^*)^n
\]

In particular, with $\sigma=0$, we recover the quantum torus $\Tscr_{h,\Ical}=[(\C^*)^d/\Z^{n-d}]$ which is dense in all these quantum toric stacks.

More generally, we can associate (see \cite{boivin2020nonsimplicial}) to a general quantum fan a stack on $\sitename{A}$ by gluing affine pieces $(\Uscr_\sigma)_{\sigma \in \Delta}$ of the form
\[
\Uscr_\sigma\coloneqq[(\C^{\dim(\sigma)} \times (\C^*)^{d-\dim(\sigma)}) \times \C^{p-d}/\Z^{n-d} \times E(\mathrm{Rel}_\C(\sigma))]
\]
where $p$ is the cardinal of a set of generators of $\sigma$ and $\mathrm{Rel}_\C(\sigma)$ is the complex vector space of relations given by these generators of $\sigma$. The quantum torus $\Tscr_{h,\Ical}$ is also dense in these stacks (cf. consequences of \cite[Definition 3.2.1.]{boivin2020nonsimplicial}).

In both cases, if $\sigma, \tau \in \Delta$ then $\Uscr_{\sigma \cap \tau}$ is an open substack of  $\Uscr_{\sigma}$ and $\Uscr_\tau$. Thereby, we get a diagram 
\begin{align}
\begin{tikzcd}[ampersand replacement=\&,scale=0.05]
    \&\& {U_\sigma} \\
	{U_{\sigma\cap\tau}} \\
	\&\& {U_\tau}
	\arrow["{\text{open }}"{pos=0.5}, hook, from=2-1, to=1-3]
	\arrow["{\text{open}}"'{pos=0.5}, hook, from=2-1, to=3-3]
\end{tikzcd}
\end{align}
for any $\sigma, \tau \in \Delta$

\begin{Def}
    The quantum toric stack associated to the fan $(\Delta,h,\Ical)$ is the colimit $\Xscr_{\Delta,h,\Ical}$ of these diagrams i.e. the gluing of the quantum toric stacks $(\Uscr_\sigma)_{\sigma \in \Delta}$ along their intersection. With stack morphisms which restrict on Picard stack morphism on tori, they form a category denoted $\mathbf{QTS}$.
\end{Def}

\begin{Thm}[\cite{katzarkov2020quantum} Theorem 6.24, \cite{boivin2020nonsimplicial} Theorem 4.2.2.2] \label{Thm : equiv_cat}
    The correspondence $(\Delta,h,\Ical) \in \mathbf{QF} \mapsto \Xscr_{\Delta,h,\Ical} \in \mathbf{QTS}$ is an equivalence of categories.
\end{Thm}
The morphism associated to a quantum fan morphism can be explicitly locally described: \\
If $(\Delta,h,\Ical)$ and $(\Delta',h',\Ical')$ are two simplicial quantum fans and $(L,H)$ is a quantum fan morphism between them, then for any cones of maximal dimension $\sigma \in \Delta, \sigma' \in \Delta'$ (with isomorphisms $L_\sigma \colon \sigma \to \Cone(e_1,\ldots,e_{d})$, $L_\sigma' \colon \sigma' \to \Cone(e'_1,\ldots,e'_{d'})$), such that $L(\sigma) \subset \sigma'$, we have a linear morphism
$A \coloneqq L_{\sigma'}LL_\sigma^{-1} \colon \R^d \to \R^{d'}$
with values in $\N^{d'}$. It defines an equivariant (for the discrete groups $\Z^{n-d}$ and $\Z^{n'-d'}$) torus morphism 
\begin{align} \label{Lbar}
    \overline{L} \colon \underline{x}\in \C^d \mapsto  \left(\underline{x}^{A^\top(e'_1)},\ldots, \underline{x}^{A^\top(e'_d)}\right) \in \C^{d'}
\end{align}
which descends to a map $\Uscr_\sigma \to \Uscr_{\sigma'}$ (the non-maximal dimension and the non-simplicial cases are a bit more involved, see \cite[Theorem 5.5]{katzarkov2020quantum} for the former and \cite[Theorem 4.1.3.3]{boivin2020nonsimplicial} for the latter)

\begin{Rem}
    If $(\Delta,h,\Ical)$ is a subfan of $(\Delta',h',\Ical')$ then $\Xscr_{\Delta,h,\Ical}$ is a open substack of $\Xscr_{\Delta',h,\Ical'}$
\end{Rem}

\begin{Rem} \label{Rem : modulispace}
Thanks to this theorem, the study of moduli spaces of quantum toric stacks can be done by the study of moduli spaces of quantum fans. In \cite{boivin2023moduli}, we consider, for $d,n \in \N$ with $n \geq d$, and $D$ a poset, the space $\Omega(d,n,D)$ of quantum fans $(\Delta,h : \R^n \to \R^d, \Ical)$ such that:
\begin{itemize}
    \item $\Delta$ ordered by inclusion is isomorphic, as poset, to $D$ ;
    \item $\Ical=[\![1,n]\!] \setminus \Delta(1)$.
\end{itemize}
We say that $\Omega(d,n,D)$ is the space of quantum fans with combinatorial type $D$ (cf. \cite[subsection 11.9]{katzarkov2020quantum}, \cite[Definition 2.1.1]{boivin2023moduli}). 
\end{Rem}

We described quantum toric stacks as the gluing of quotient stacks. Fortunately, they can also be described as global quotient stacks.

\begin{Thm}[\cite{katzarkov2020quantum} Theorem 7.6, Quantum GIT] \label{Thm : QGIT}
If $(\Delta,h,\Ical)$ is a simplicial quantum fan then the quantum toric stack $\Xscr_{\Delta,h,\Ical}$ is a quotient stack
\[
\Xscr_{\Delta,h,\Ical}=[\Sscr(\Delta)/\C^{n-d}]
\]
where $\Sscr(\Delta)$ is a quasi-affine (classical) toric variety given by the combinatorics of $\Delta$ :
    \begin{equation} \label{Def_S}
        \Sscr(\Delta)=\bigcup_{\sigma \in \Delta_{max}} \C^{\sigma(1)} \times (\C^*)^{\sigma(1)^c} \subset \C^n;
    \end{equation}    
and $\C^{n-d}$ acts on $\Sscr$ through 
\[
t \cdot z\coloneqq E(k \otimes id_\C(t))z
\]
where $k$ is a Gale transform of $h \otimes id_\R : \R^n \to \R^d$ i.e. a linear map such that the short sequence 
\begin{align}\label{Gale_transform}
\begin{tikzcd}[ampersand replacement=\&]
	0 \& {\R^{n-d}} \& {\R^n} \& {\R^d} \& 0
	\arrow[from=1-1, to=1-2]
	\arrow["k", from=1-2, to=1-3]
	\arrow["h", from=1-3, to=1-4]
	\arrow[from=1-4, to=1-5]
\end{tikzcd}
\end{align}
is exact
\end{Thm}

\begin{Rem}
    Let $(L,H) \colon (\Delta,h,\Ical) \to (\Delta',h',\Ical')$ be a quantum fan morphism. Then $H$ is a fan morphism $\Delta_h \to \Delta_{h'}$ where $\Delta_h$ is the fan of $\R^n$ such that $\Cone(e_1,\ldots,e_k) \in \Delta_h$ if $\Cone(h(e_1),\ldots,h(e_k)) \in \Delta$. Hence, it defines a toric morphism
    \begin{align} \label{Hbar}
    \overline{H} \colon \Sscr(\Delta)=X_{\Delta_h} \to X_{\Delta'_{h'}}=\Sscr(\Delta').
    \end{align}
\end{Rem}

\subsection{Fibers of a toric morphism}
\label{subsection : fibers}
In this subsection, we compute the fiber over a torus-equivariant substack of a toric morphism between quantum toric stacks. It will be useful later when we study the exceptional divisor of a birational morphism.

\begin{Thm} \label{Thm : cal_fiber}
Let $\lcal : \Uscr_\sigma\coloneqq [\C^d/\Z^{n-d}] \to \Uscr_\tau \coloneqq [\C^{d'}/\Z^{n'-d'}]$ be a toric morphism induced by the linear morphisms $(L,H)$. Then, 
\[\lcal^{-1}(0)\coloneqq [\{0_{\C^{n'-d'}}\}/ \Z^{n'-d'}] \times_{[\C^{d'}/\Z^{n'-d'}]} [\C^{d}/\Z^{N-d}]=[\overline{L}^{-1}(0)/\Z^{n-d}]
\]
where $\overline{L}$ is the map $\C^d \to \C^d$ induced by $L$ described in \eqref{Lbar}. In other words, the following diagram is cartesian: 
\[\begin{tikzcd}
	{[\overline{L}^{-1}(0)/\Z^{n-d}]} & {[\C^d/\Z^{n-d}]} \\
	{[0/\Z^{n'-d'}]} & {[\C^{d'}/\Z^{n'-d'}]}
	\arrow[shift right=1, hook, from=1-1, to=1-2]
	\arrow["\lcal", from=1-2, to=2-2]
	\arrow[from=1-1, to=2-1]
	\arrow["\lrcorner"{anchor=center, pos=0.125}, draw=none, from=1-1, to=2-2]
	\arrow["i",shift right=1, hook, from=2-1, to=2-2]
\end{tikzcd}\]
\end{Thm}

\begin{proof}
Let $T \in \sitename{A}$. Then,
\[
\lcal^{-1}(0)(T)=\left\{ \alpha : T' \stackrel{\Z^{n'-d'}}{\longrightarrow} T,  \begin{tikzcd}
	{\widetilde{T}} & {\C^d} \\
	T
	\arrow["\pi"', from=1-1, to=2-1]
	\arrow["m", from=1-1, to=1-2]
\end{tikzcd} , \varphi \colon i(\alpha) \simeq \lcal(\pi,m) \right\}
\]
Set $\widehat{T}\coloneqq \widetilde{T}\times_{H_2} \Z^{n-d}$. \\
The isomorphism $\varphi$ (in $[\C^{d'}/\Z^{N'-d'}](T)$) is described by the following commutative diagram:

\[\begin{tikzcd}
	& {\C^{d'}} \\
	{T'} && {\widehat{T}} \\
	& T
	\arrow[from=2-1, to=3-2]
	\arrow[from=2-3, to=3-2]
	\arrow["\varphi", from=2-1, to=2-3]
	\arrow["0", from=2-1, to=1-2]
	\arrow[from=2-3, to=1-2]
	\arrow["\simeq"', from=2-1, to=2-3]
\end{tikzcd}\]
where the map $m' \colon \widehat{T}\to \C^{d'}$ is defined by:
\[
m'([\widetilde{t},p])=\overline{L}(m(\widetilde{t}))E'(h(0\oplus p))
\]
Since the upper triangle is commutative then, for all $\widetilde{t} \in \widetilde{T}'$, 
\[
\overline{L}(m(\widetilde{t}))=0
\]
In other words, $m$ takes values in $\overline{L}^{-1}(0)$.  

\begin{Lemme} \label{lemme : forgetful_equiv}
For all object $T$ in $\sitename{A}$, the forgetful functor $\lcal^{-1}(0)(T) \to [\overline{L}^{-1}(0)/\Z^{n-d}](T)$ is an equivalence of categories.
\end{Lemme}
\begin{proof}
This functor is clearly (essentially) surjective. It remains to us to prove that it is fully faithful: \\
A morphism $(\alpha,(\pi,m),\varphi) \to (\alpha',(\pi',m'),\varphi')$ of $\lcal^{-1}(0)(T)$ is given by a morphism $u : \alpha \to \alpha'$ and a morphism $v : (\pi,m) \to (\pi',m')$ such that the following diagram commutes: 
\begin{equation} \label{morphism_fiber}
\begin{tikzcd}
	{i(\alpha)} && {i(\alpha')} \\
	{\lcal(\pi,m)} && {\lcal(\pi',m')}
	\arrow["{i(u)}", from=1-1, to=1-3]
	\arrow["{\lcal(v)}", from=2-1, to=2-3]
	\arrow["\varphi"', from=1-1, to=2-1]
	\arrow["{\varphi'}", from=1-3, to=2-3]
\end{tikzcd}
\end{equation}
In order to prove that this forgetful functor is fully faithful, it remains to us to prove that the map
\[ F(\--)\colon (u,v)\in \Hom((\alpha,(\pi,m),\varphi),(\alpha',(\pi',m'),\varphi')) \mapsto v \in \Hom((\pi,m),(\pi',m')), 
\]
is bijective.\\
In the diagram \eqref{morphism_fiber}, we can remark that $i(u)=(\varphi')^{-1}\lcal(v) \varphi$. By corestricting the equivariant map of $i(u)$ on 0, we get an element, denoted $u(v)$, of $[0/\C^{n'-d'}](T)$. \\It defines a map 
\[v \in \Hom((\pi,m),(\pi',m'))\mapsto \left(u(v), v \varphi\right)\in\Hom((\alpha,(\pi,m),\varphi),(\alpha',(\pi',m'),\varphi'))\] which is the inverse of $F(\--)$.
\end{proof}
The forgetful functor defines a stack morphism $\lcal^{-1}(0) \to [\overline{L}^{-1}(0)/\Z^{n-d}]$ which is an isomorphism thanks to lemma \ref{lemme : forgetful_equiv}.
\end{proof}

\begin{Cor}
We can replace $\{0\}$  by any invariant subvariety of $\C^n$ in the statement of theorem \ref{Thm : cal_fiber} and the $\C^d$ by $\C^k \times (\C^*)^{d-k}$, $0 \leq k \leq d$.
\end{Cor}

\begin{Not}
If $t=(t_1,\ldots,t_d) \in \C^d$ (or\footnote{This is the same thing since the action preserves the zero coordinates}
if $[t]$ is the orbit of an element of $\C^d$), we note $I_t\coloneqq\{ i \in \{1,\ldots, d\} \mid t_i \neq 0\}$ and $t_{\neq 0}=(t_i)_{i \in I_{t}}$.
\end{Not}

We continue by computing the inverse image of an orbit: 

\begin{Lemme} \label{Lemme : fiber_computation}
Let $\lcal$ be a toric morphism $[\C^d/\Z^{n-d}] \to [\C^{d}/\Z^{n'-d}]$. Let $L=(\alpha_{j}^i)_{1 \leq i,j \leq d} \colon \R^d \to \R^d, H : \Z^n \to \Z^{n'}$ be the linear morphisms associated to $\lcal$. Let $t \in \C^d$. 

Then the reduced analytic space $X \coloneqq \overline{L}^{-1}([t])^{\mathrm{red}}$ is 
\begin{itemize}
    \item $E(L_\C^{-1}(w))$ if $t=E(w) \in (\C^*)^d$ .
    \item $\bigcup_A \{0\}^A \times E(L_A^{-1}(w))$ where the union is indexed over all the subsets $A$ of $\{1,\ldots,d\}$ such that
    \begin{itemize}
        \item $\forall i \in I, \forall l \in A, \alpha_l^i=0$
        \item $\forall k \in I^c, \exists j \in A, \alpha_j^k\neq 0$
    \end{itemize}
    where $L_A\coloneqq (\alpha_j^i)_{i \in I_{[t]}, j \in A^c}$ if $t=(0_{I_{[t]}},E(w))$.
    \item $\bigcup_A \C^{A^c} \times 0_{A} $ where the union is indexed by the subsets of $\{1,\ldots,d\}$ such that for all $i \in [\![1,d]\!]$, there exists $j \in A$, $\alpha_j^i \neq 0$, if $t=0$
\end{itemize}
\end{Lemme}

\begin{proof}
If $t \in (\C^*)^d$ then, since the restriction of a toric morphism on tori is a torus morphism, we obtain our result using the commutative diagram \eqref{induced_torus_morph}.\\
If $t \in \C^d \setminus (\{0_d\} \cup (\C^*)^d)$ then $I_t \neq \emptyset$ and $I_t \neq \{1,\ldots,d\}$. If $z=(z_1,\ldots,z_n) \in X$ then the coordinates of $z$ verify the following conditions:
\begin{enumerate}
    \item For all $i \in I_{[t]}$ and all $l \in \{1,\ldots,d\}$, $\alpha_l^i=0$ or $z_l \neq 0$ (the elements of non-zero product is non-zero)
    \item For all $k \notin I_{[t]}$, there exists $j \in \{1,\ldots,d\}$ such that $\alpha^k_{j} \neq 0$ et $z_j=0$ (a product is zero if, and only if, one element is zero).
\end{enumerate}
If we note $A$ the set of zero coordinates of $z$, we get the following conditions :
\begin{enumerate}
    \item For all $i \in I_{[t]}$ and all $l \in A$, $\alpha_l^i=0$ 
    \item\label{condition_2} For all $k \notin I_{[t]}$, there exists $j \in A$ such that $\alpha^k_{j} \neq 0$
\end{enumerate}

Then, for the remaining coordinates, we conclude in the same way as in the first point.\\
If $t=0$ then $z \in X$ if, and only if, the condition \eqref{condition_2} is verified. Then, for all subset $A$ of $\{1,\ldots,d \}$ such that, for all $i \in [\![1,d]\!]$, there exists $j \in A$, $\alpha_j^i \neq 0$, $(z_{A^c},0_{A}) \in X$ for any $z_{A^c} \in \C^{A^c}$. Conversely, every point of $X$ are of this form.

\end{proof}

\begin{Ex}
Consider $\left(L=\begin{pmatrix} a & b \\c & d\end{pmatrix} \colon \C^2 \to \C^2,H : \Z^n \to \Z^{n'}\right)$ (where $a,b,c,d\in \N_{\geq 0}$ and  $ad-bc \neq 0$) be a pair of linear maps compatible with two calibrations $\Z^n \to \Gamma$ and $\Z^{n'} \to \Gamma'$ in $\R^2$. The map $\overline{L}$ is given by:
\[
(z_1,z_2) \mapsto (z_1^az_2^b,z_1^cz_2^d)
\]
and descend to quotient in a morphism $\lcal : [\C^2/\Z^{n-2}] \to [\C^2/\Z^{n'-2}]$.

Then, 
\begin{itemize}
    \item If $t=\Z^{n-2} \cdot (E(w_1),E(w_2))  \in (\C^*)^2/\Z^{n'-2}$ then thanks to the following diagram
    \[\begin{tikzcd}
	{\C^2} & {\C^2} \\
	{[(\C^*)^2/\Z^{n-2}]} & {[(\C^*)^2/\Z^{n'-2}]}
	\arrow[from=1-1, to=2-1]
	\arrow["L", from=1-1, to=1-2]
	\arrow[from=1-2, to=2-2]
	\arrow["{\lcal}", from=2-1, to=2-2]
\end{tikzcd}\]
we see that the fiber over $[t]$ is $[E(L^{-1}(w))]$ (the compatibility with the quotients comes from the compatibility with the linear maps) ;
\item if $t=\Z^{n-2}\cdot (0,E(w))$ then there are several possibilities: 
\begin{itemize}
    \item If $a=0$ then $b \neq 0$ (and if $b=0$ then $a \neq 0$) and $z_2=0$ (resp. $z_1=0$) then $d=0$ (resp. $c=0$) otherwise we would have a contradiction with the non-nullity of the second (resp. first) coordinate. All that remains is to solve $[z_1^c]=[E(w)]$ (resp. $[z_2^d]=[E(w)]$). The reduced fiber of $L$ over $z$ is therefore $[\mu_c E(w/c)] \times 0$ (resp. $0 \times[\mu_d E(w/d)] $) where $\mu_n$ is the group of $n$th roots of unity;
    \item if $a \neq 0$ and $b \neq 0$ then $z_1=0$ or $z_2=0$ and we arrive at the same conclusion as in the previous case;
\end{itemize}
\item if $t=(0,0)$ then the same reasoning applies as above, but $E(t)$ is replaced by $0$, so the reduced fiber is $0$ if the power is non-zero.
\end{itemize}

\end{Ex}

\section{Toy case: blow-up of a quantum plane at the origin}
\label{section : toycase}

Let $(\Z^3 \to \Gamma \subset \R^2,\emptyset)$ be a (standard) calibration such that $v=(a,b)\coloneqq h(e_3) \in \R_{> 0}^2$.

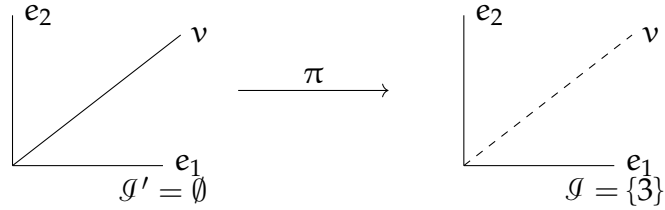
\begin{figure}[!ht]
    \centering
    \begin{tikzpicture}
    \draw (0,0) to (0,2) node[right]{$e_2$} ;
    \draw (0,0) to (2,0) node[right]{$e_1$}
    node[left, below]{$\Ical'=\emptyset$} ;
    \draw (0,0) to (2.236,1.738)node[right]{$v$}   ;
    \draw[->] (3, 1) to (5,1) node at (4,1.2) {$\pi$} ; 
    \draw (6,0) to (6,2) node[right]{$e_2$} ;
    \draw (6,0) to (8,0) node[right]{$e_1$} 
     node[left,below]{$\Ical=\{3\}$} ;
    \draw[dashed](6,0) to (6+2.236,1.738)node[right]{$v$}   ;
\end{tikzpicture} 
\caption{blow up of a plane at the origin}
    \label{fig:eclat_plan}
\end{figure}

Classically, a blow-up morphism is described by the identity morphism between two compatible fans. The natural blow-up morphism would be described as in figure \ref{fig:eclat_plan} and induced by the pair $(id_{\R^2},H=\begin{pmatrix}
1 & 0 & a \\
0 & 1 & b \\
0 & 0& 0 
\end{pmatrix} : \R^3 \to \R^3
)$

The additional conditions on quantum fan morphisms lead to the following statement:
\begin{Lemme}
The "natural blow-up morphism" $\pi$  is a quantum fan morphism
if, and only if, $v \in \N^2$. 
\end{Lemme}
\begin{proof}

This comes from the definition of fan morphisms and the equivalence of categories of theorem \ref{Thm : equiv_cat}. 
\end{proof}

\begin{Rem} \label{Rem : non_iso}
Note that the induced morphism by $\pi$ is never an isomorphism (since $H \notin \GL_3(\Z)$) even if we restrict on (equivariant) open substacks.
\end{Rem}

However, we can prove that the morphism is birational as stack morphisms: 
\\
The morphism $\pi$ induces a morphism $\overline{H} : \Sscr(\Delta')=\C^2 \setminus \{0\} \times \C \to \Sscr(\Delta)=\C^2 \times \C^* $ (see equation \eqref{Hbar})  defined by:
\[
\overline{H}(x,y,z)=(xz^a,yz^b,1)
\]
The quantum toric stack on the left (resp. right) is the quotient of $\C^2 \setminus \{0\} \times \C$ (resp. of $\C^2\times \C^*$) by the action of $\C$ given by
\[
t \cdot (z_1,z_2,z_3)=(E(-at)z_1,E(-bt)z_2,E(t)z_3)
\]

We conclude that for any $(x,y,z) \in \C^2 \setminus \{0\} \times \C^*$, $H(x,y,z)$ is in the same equivalence class as $(x,y,z)$ (indeed, if $z=E(t)$ then $(xz^a,yz^b,1)=(-t) \cdot (x,y,z)$). Consequently, the blow-up morphism induced by $\pi$ can be restricted to the identity on $[\C^2 \setminus \{0\} \times \C^*/\C]$ (hence it is a birational morphism but the inverse is not naturally a toric morphism). Note that the associated morphism $\R^3 \to \R^3$ is also the identity. Hence, in contrast with \ref{Rem : non_iso}, we have a quantum fan isomorphism.

\begin{Rem} \label{Lemme : eclatement2D}
The identity morphism is always well-defined for any $v \in\R_{> 0}^2$. However, it is not a morphism $\Sscr(\Delta') \to \Sscr(\Delta)$ but only on open subsets i.e. $\C^2 \setminus \{0\} \times \C^* \subset \Sscr(\Delta') \to (\C^2 \times \{0\}) \times \C^* \subset \Sscr(\Delta)$
\end{Rem}

\section{Rational weighted blow-up}
\label{section : rationalwbu}
\subsection{Integer weight case}

\label{integer_weight}

In this paragraph, we will consider the case of the integer-weight blow-up\footnote{
Note that the "classical" weighted blow-up introduced in \cite{Abramovich_2010} (see also \cite[section 10.2.7]{olsson2016algebraic}, and \cite[Example 5.4.2]{w_blowup} for the toric stacky case) is a $\Z^{n-\mathrm{rk}(\Gamma)}$ gerbe over this blow-up.} of an equivariant closed analytic substack of an affine quantum toric stack (given by a simplicial cone), in other words the blow-up of a substack $[\C^I\times 0^{I^c}/\Z^{n-d}] $ in $\Uscr \coloneqq [\C^d/\Z^{n-d}]$:

Let $I$ be a subset of $\{1,\ldots,d\}$ such that $\Card(I) \leq d-2$. Let $v=\sum_{k \in I} v_k e_k$ be a vector of $(\N_{>0})^{I^c} \times 0^I \subset C_d=\R_{\geq 0}^d$ such that there exists $k$ such that $h(e_k)=v$. Up to reordering, we can assume that $k=d+1$.
For all $i \in I^c$, let $\sigma_i$ be the cone of $\R^d$ defined by 
\[
\sigma_i=\Cone(e_1,\ldots,\widehat{e_i},v,e_{i+1},\ldots,e_d).
\]
Note $\Delta^*(v)$ the fan whose maximum cones are the $\sigma_i$. 
Note $\fcal$ the toric morphism given by the quantum fan morphism $\Delta^*(v) \to C_d$ given by the pair $(id,H)$ where
\[
H\coloneqq \begin{pmatrix}
1&     &0 &v_1      &  0  &\cdots&0\\
&\ddots& & \vdots  &  \vdots  && \vdots\\
0&      &1& v_d     &  0  &\cdots &0\\
0&   \cdots   & \cdots  &    0     &   1&&0 \\
\vdots&      & &    \vdots     &   &\ddots& \\
0&  \cdots     &\cdots  &    0     & 0  &&1 \\
         
\end{pmatrix}.
\]

\begin{Rem}
It is a generalization of the example of section \ref{section : toycase}. And in the same way as \ref{Rem : non_iso}, the matrix $H$ is not invertible.
\end{Rem}

\begin{Def} 
    The morphism $\fcal\colon \Xscr_{\Delta^*(v),h,\Ical \setminus \{d+1\}} \to \Uscr$ is called blow-up morphism of $\Uscr$ at $[0^I \times \C^{I^c}/\Z^{n-d}]$. 
\end{Def}

We now describe locally the fiber of this morphism. We postponed the global description to the section \ref{general_blowup}.

Let $i \in I^c$. Let $L_i : \R^d \to \R^d$ be the linear isomorphism such that $L_i(e_k)=e_k$ for $k \neq i$ and $L_i(v)=e_i$ and $H_i$ the permutation $(d+1 \; i) \in \Sfrak_n$ . \\
Let $\Uscr_i=[\C^d/\Z^{n-d}]$ (where $\Z^{n-d}$ acts on $\C^d$ via $E L_i h H_i^{-1}$) be the quantum toric stack associated with $\sigma_i$. \\
The quantum fan morphism \[C_d \xrightarrow[\simeq]{(L_i^{-1},H_i^{-1})}  \sigma_i \xhookrightarrow[]{(id,H)} C_d\] induces a toric morphism $ \fcal_i : \Uscr_i\to \Uscr$.

Note $A_i=\begin{pmatrix}
1 & & &v_1 & 0& \cdots& 0\\
 & \ddots && \vdots &\vdots& & \vdots\\
 & &1& \vdots &\vdots && \vdots\\
0 & \cdots &0& v_i & 0& \cdots&0\\
\vdots & & \vdots& \vdots & 1 & &  0\\
 \vdots & &\vdots& \vdots & & \ddots & \\
0 & \cdots&0& v_d & 0&& 1 \\\\
\end{pmatrix} : \C^d \to \C^d$ the inverse of $L_i$. The morphism $A_i$ descends into a morphism $\overline{A_i}:(\C^*)^d=\C^d/\Z^d \to (\C^*)^d=\C^d/\Z^d$ which can be extended into a toric morphism $\overline{A_i} : \C^d \to \C^d$ defined by : 
\[\forall z=(z_1,\ldots, z_d) \in \C^d, \overline{A_i}(z)=(z_1 z_i^{v_1},\ldots,z_{i-1}z_i^{v_{i-1}},z_i^{v_i},z_{i+1}z_i^{v_{i+1}},\ldots,z_{d}z_i^{v_{d}})
\]

Then the fiber $\fcal_i^{-1}([a])$ is the quotient of $\overline{A_i}^{-1}([a])$. For all $i,j \in I$, the fibers $\fcal_i^{-1}([a])$ and $\fcal_j^{-1}([a])$ coincide on the open substack $\Uscr_{\sigma_i \cap \sigma_j}$. The fibers $(\fcal_i^{-1}([a]))_{i \in I}$ glue into the fiber over $[a]$ of $\fcal$, which will be a quotient of $\overline{H}^{-1}([a])$ thanks to the theorem \ref{Thm : QGIT}.

\begin{Lemme}
The fiber of the morphism $\overline{A}_i$ over a point
$\alpha \in 0^{I^c} \times \C^I$ is the analytic space
\[
\overline{A_i}^{-1}(\alpha)=V(z_j-\alpha_j, j \in I,z_i^{v_i},z_kz_i^{v_k}, k \in I^c \setminus \{i\}) 
\]
\end{Lemme}

The projection
\[V(z_j-\alpha_j, j \in I,z_i^{v_i},z_kz_i^{v_k}, k \in I^c \setminus \{i\}) \to V(z_i^{v_i},z_kz_i^{v_k}, k \in I^c \setminus \{i\}) \eqqcolon V\]
is an isomorphism. \\
By lemma \ref{Lemme : fiber_computation}, we have :

\begin{Lemme} \label{Lemma : underlyingfibrespace}
The underlying topological space of $V$ is $\C^{I^c \setminus \{i\}} \times 0^{\{i\}} $ 
\end{Lemme}

Thanks to its explicit equations, we can describe the stalks of the structure sheaf of $V$:

\begin{Lemme} \label{Lemme : calcul_fibre}
The stalk of the structure sheaf $\Ocal$ of $V$ on the point $\beta=(\beta_j, j \in I^c \setminus \{i\},0)$ is reduced if, and only if, $\min\limits_{j \in J} v_j=1$ and if for all $k \in I^c \setminus {J}$, $v_k=1$ where $J$ is the set of indices of $I^c\setminus \{i\}$ of non-zero coordinates of $\beta$. In particular, we have: 
\begin{itemize}
    \item If $\beta \in (\C^*)^{I^c \setminus \{i\}}$ then $\Ocal_\beta$ is reduced if, and only if, $\min\limits_{j \in I^c} v_j=1 $ ; 
    \item If $\beta=0$ then $\Ocal_0$ is reduced if, and only if, $v=\sum_{k \in I^c \setminus \{i\}} e_k$.
\end{itemize}

\end{Lemme}

The previous statements prove

\begin{Thm}
The analytic space $\overline{A_i}^{-1}(\alpha)$ is reduced if, and only if,  $v=\sum_{k \in I^c \setminus \{i\}} e_k$.
\end{Thm}

By gluing the fibers of the $\overline{A_i}^{-1}(\alpha)$, we obtain the following statement: 

\begin{Thm}
The blow-up morphism $\fcal : \Xscr_{\Delta^*(v),h,\Ical} \to \Uscr$ has reduced fibers if, and only if, $v=\sum_{k \in I^c} e_k$
\end{Thm}

\subsection{Rational weight case}

In this subsection, we keep the notation of the subsection \ref{integer_weight} except that we assume now $v \in \Q_{>0}^{I^c} \times \{0\}^I$. The pair of maps $(id,H)$ is no longer a fan morphism but induces a birational stack morphism (but has not a quantum toric inverse) cf. section \ref{section : toycase}) \\
There exists $N \in \N_{>0}$ such that $Nv \in (\N_{>0})^{I^c} \times 0^I$. Then the cones of $N \widetilde{C_d}$ are generated by elements of $\N^d$. Then, since $N\Gamma$ is a subgroup of $\Gamma$, the identity morphism of $\R^d$ and the morphism 
\[
NH=\begin{pmatrix}
N & & 0&Nv_1 &0 & \cdots &0\\
 & \ddots && \vdots & \vdots & & \vdots \\
0 &  &N& Nv_d & \vdots & & \vdots\\
0 & \cdots&0& 0 &0 & \cdots & 0\\
\vdots & &\vdots& \vdots & N & & 0\\
\vdots & &\vdots& \vdots & & \ddots & \\
0 & \cdots&0& 0 & 0& & N \\
\end{pmatrix}
\]
induced a toric morphism
\[
\Xscr_{N h,N \widetilde{C_d}} \longrightarrow \Uscr_{h,C_d}
\]
We thus obtain a first example of the decomposition of the birational morphism into a zig-zag of toric morphisms: 

\[\begin{tikzcd}[ampersand replacement=\&]
	{\Xscr_{h,\widetilde{C_d}}} \&\& {\Uscr_{C_d}} \\
	\& {\Xscr_{Nh,N\widetilde{C_d}}}
	\arrow["{(id,H)}", dashed, from=1-1, to=1-3]
	\arrow["{(id,NH)}"', from=2-2, to=1-3]
	\arrow["{(id,Nid)}", from=2-2, to=1-1]
\end{tikzcd}\]

The toric morphism associated with the pair $(id,NH)$ is given by a family of toric morphisms $\lcal_i : \Uscr_{N\sigma_i} \to \Uscr_{h,C_d}$ and thus, by linear maps: 
\[
L_i=\begin{pmatrix}
N & & 0&Nv_1 &0 & \cdots&0 \\
 & \ddots && \vdots &\vdots & &\vdots \\
0 &  &N& \vdots & \vdots& &\vdots \\
0 & \cdots&0& Nv_i &0 & & 0\\
\vdots & &\vdots& \vdots & N & &0 \\
\vdots & &\vdots& \vdots & & \ddots & \\
0 & \cdots &0& Nv_d &0 & & N \\
\end{pmatrix}
\]
Each of these toric morphisms is given by the descent to the quotient of toric morphisms $\overline{L}_i \colon \C^d \to \C^d$ defined by
\[
\forall z=(z_1,\ldots,z_d) \in \C^d, \overline{L}_i(z)=(z_1^N z_i^{Nv_1},\ldots,z_{i-1}^Nz_i^{Nv_{i-1}},z_i^{Nv_i},z_{i+1}^N z_i^{N v_{i+1}},\ldots,z_{d}z_i^{N v_{d}})
\]
\begin{Lemme}
If $v \notin \N^{I^c}$ then for all $\alpha \in 0^{I^c} \times \C^I$, the analytic space $\overline{L_i}^{-1}(\alpha)$ is not reduced.
\end{Lemme}
\begin{proof}
In this case, $N>1$ and the statement follows from lemma \ref{Lemme : calcul_fibre}.
\end{proof}

\section{General weighted blow-up}
\label{general_blowup}
\subsection{Definition and properties of birational fan morphisms}
 The example in the previous section illustrates that in order to define blow-ups, it is necessary to consider linear maps $(L,H)$ that are fan morphisms on sub-fans: 
 
 \begin{Def} \label{Def : morph_birat}

A birational quantum fan morphism between $(\Delta,h \colon \R^n \to \R^d,\Ical)$ in $\Gamma \subset \R^d$ and $(\Delta',h'\colon \R^n \to \R^d,\Ical')$ in $\Gamma'$ is a pair of linear morphisms $(L: \R^d \to \R^{d},H: \R^n \to \R^{N})$ such that there exists a pair of subfans $\widetilde{\Delta} \subset \Delta$, $\widetilde{\Delta'} \subset \Delta'$, subsets $J \subset [\![1,n]\!] \setminus (\Delta(1) \cup \Ical)$, $J' \subset [\![1,n']\!] \setminus (\Delta'(1) \cup \Ical')$ such that $(L,H)$ is a quantum fan isomorphism between $(\widetilde{\Delta},h,\Ical \cup (\Delta(1) \setminus \widetilde{\Delta}(1)) \cup J)$ and $(\widetilde{\Delta'},h,\Ical' \cup (\Delta(1) \setminus \widetilde{\Delta}'(1)) \cup J')$. \\
We note this morphism $(\Delta,h,\Ical)\dashrightarrow (\Delta',h',\Ical')$.
\end{Def}

\begin{Rem}
    In particular, the morphisms $L$ and $H$ are isomorphisms.
\end{Rem}

\begin{Rem}
We allow change of generators which are neither fan generators nor virtual generators, so that the morphisms induced by the cobordisms of section \ref{chap:10} are birational morphisms.
\end{Rem}

The symmetry of this definition gives us the following result:
\begin{Lemme}
If $(L,H)$ is a birational quantum fan morphism $(\Delta,h,\Ical)\dashrightarrow (\Delta',h',\Ical')$ then $(L^{-1},H^{-1})$ is a birational quantum fan morphism $(\Delta',h',\Ical')\dashrightarrow (\Delta,h,\Ical)$.
\end{Lemme}

By the equivalence of categories described in theorem \ref{Thm : equiv_cat}, a birational fan morphism $(\Delta,h,\Ical)\dashrightarrow (\Delta',h',\Ical')$ induces a birational morphism $\Xscr_{\Delta,h,\Ical}\dashrightarrow \Xscr_{\Delta',h',\Ical'}$ i.e. a toric isomorphism between two open substacks $\Uscr \subset \Xscr_{\Delta,h,\Ical}$ and $\Vscr \subset \Xscr_{\Delta',h', \Ical'}$ which are given by the fans $(\widetilde{\Delta},h,\Ical \cup (\Delta(1) \setminus \widetilde{\Delta}(1)) \cup J)$ and $(\widetilde{\Delta'},h,\Ical' \cup (\Delta'(1) \setminus \widetilde{\Delta}'(1)) \cup J')$.

\begin{Ex} 
Let $\alpha=(\alpha_1,\alpha_2) \in \R_{>0}^2 \setminus \Q^2$ and $h \colon (x,y,z)  \in \R^3 \mapsto (x+\alpha_1 z,y+\alpha_2 z) \in \R^2$. The $\alpha$-blow-up of the plane at the origin is defined by the pair 
\[(id_{\R^2}, id_{\R^3}) \colon (\widetilde{\Delta},h,\emptyset) \dashrightarrow (\Delta_{\R_{\geq 0}^2},h,\{3\})\]
(where $\Delta_{\R_{\geq 0}^2}$ is the fan of faces of the cone $\R_{\geq 0}^2$ and $\widetilde{\Delta}$ is the fan whose maximal cones are $\Cone(e_1,\alpha)$ and $\Cone(\alpha,e_2)$, see figure \ref{fig:eclat_plan} for an illustration)  which induces a birational morphism $\Xscr'=[\C^2 \setminus \{0\} \times \C/\C] \dashrightarrow \Xscr=[\C^2 \times \C^*/\C]$ where the action of $\C$ is given by 
\[ 
t \cdot (z_1,z_2,z_3)=(E(-\alpha_1t)z_1,E(-\alpha_2 t)z_2,E(t)z_3)
\]
As seen in lemma \ref{Lemme : eclatement2D}, the toric open substacks in question are $\Uscr=\Vscr=[\C^2 \setminus \{0\} \times \C^*/\C]$ and the morphism between them is the identity. The complementary $\Xscr' \setminus [\C^2 \setminus \{0\} \times \C^*/\C]=[\C^2 \setminus \{0\} \times \{0\}/\C]$ is the "exceptional divisor" of the blow-up (it is a quantum projective line) and $\Xscr \setminus [\C^2 \setminus \{0\} \times \C^*/\C]=[\{0\} \times \C^*/\C]\simeq \Bscr \Z=[*/\Z]$
is the "point" we have blown-up.

\end{Ex}

\begin{Prop}
The composition of two birational quantum fan morphisms is a birational quantum fan morphism.
\end{Prop}
\begin{proof}
Let $(L_1,H_1) : (\Delta,h,\Ical) \dashrightarrow (\Delta',h',\Ical')$ and $(L_2,H_2) : (\Delta',h',\Ical') \dashrightarrow (\Delta'',h'',\Ical'')$ be two birational fan morphisms. They induce quantum fan isomorphisms \[(\widetilde{\Delta},h,\Ical \cup (\Delta(1) \setminus \widetilde{\Delta}(1)) \cup J) \simeq (\widetilde{\Delta'},h,\Ical' \cup (\Delta(1) \setminus \widetilde{\Delta}'(1)) \cup J')\] and 
\[(\widetilde{\widetilde{\Delta}}',h',\Ical' \cup (\Delta'(1) \setminus \widetilde{\widetilde{\Delta'}}(1)) \cup J) \simeq (\widetilde{\Delta''},h,\Ical'' \cup (\Delta''(1) \setminus \widetilde{\Delta}''(1)) \cup J').\]
Let $(\overline{\Delta}',h,\overline{\Ical}')$ be a common refinement of the fans $(\widetilde{\Delta'},h,\Ical' \cup (\Delta(1) \setminus \widetilde{\Delta}'(1)) \cup J')$ and $(\widetilde{\widetilde{\Delta}}',h',\Ical' \cup (\Delta'(1) \setminus \widetilde{\widetilde{\Delta'}}(1)) \cup J)$. With the isomorphisms $(L_1^{-1},H_1^{-1})$ and $(L_2,H_2)$, we get compatible sub-fans of $(\Delta,h,\Ical)$ and $(\Delta'',h'',\Ical'')$. We conclude that $(L_2L_1,H_2H_1)$ is a birational morphism of inverse $(L_1^{-1}L_2^{-1},H_1^{-1}H_2^{-1})$.
\end{proof}
\subsection{Irrational weighted blow-up}
\label{blow_up_section}

This subsection is devoted the description of the irrational-weighted blow-up of a $\Tscr_{h,\Ical}$-equivariant substack of $\Uscr\coloneqq[\C^d/\Z^{n-d}]$.

The starting point is almost the same as in subsection \ref{integer_weight} : \\
Let $d \in \N_{>0}$ be a natural number and $\Gamma$ be a finitely generated subgroup of $\R^d$. Let $(h : \Z^n \to \Gamma,\Ical)$ be a standard calibration of $\Gamma$.
Let $I$ be a subset of $\{1,\ldots,d\}$ such that $\Card(I) \leq d-2$. Let $i_1^c < \ldots < i^c_{\Card(I^c)}$ be the elements of $I^c$. Let $\alpha=\sum_{j \in I^c} \alpha_j e_j$ a vector of $\R^{I^c}_{>0} \times 0^I \subset C_d$ such that there exists $k$ such that $h(e_k)=\alpha$. \\
For any $1 \leq j \leq \Card(I^c)$, let $\sigma_j$ be the cone of $\R^d$ defined by 
\[
\sigma_j\coloneqq \Cone(e_k, k \in I, e_{i^c_1},\ldots, \widehat{e_{i^c_j}},\alpha, e_{i^c_{j+1}},\ldots, e_{i^c_{\Card(I^c)}})
\]
Note $(\widetilde{C_d},h,\Ical) $ be the fan whose maximum cones are the $\sigma_j$.
\\
Let $i \in I^c$. Let $A_i : \R^d \to \R^d$ the linear isomorphism such that $A_i(e_k)\coloneqq e_k$ for $k \neq i$ and $A_i(e_i)\coloneqq\alpha$ and $H_i$ the permutation $(k \; i) \in \Sfrak_n \subset \GL_n(\R)$. 
\\
Let $\Uscr_i=[\C^d/\Z^{n-d}]$ (where $\Z^{n-d}$ acts on $\C^d$ through $E A_i^{-1} h H_i$) be the quantum toric stack associated with $\sigma_i$.
\\
The birational fan morphism $(id_{\R^d},id_{\R^n})$ (see \ref{Lemme : eclatement2D} for this choice of morphisms) induces a birational morphism ---called blow-up of weight $\alpha$---
\[
\fcal \colon \Xscr_{\widetilde{C_d},h,\Ical} \dashrightarrow \Uscr_{C_d}.
\]

Its restriction on the complementary of the exceptional divisor induces an isomorphism between
\[
\Xscr_{\widetilde{C_d},h,\Ical} \setminus \Dscr \simeq [\C^d \setminus \{0\}^{I^c} \times \C^I/\Z^{n-d}].
\]
The birational inverse of $\fcal$ is called blow-down (of weight $\alpha$).

The morphism $\fcal$ induces a family $\fcal_i : \Uscr_{i} \to \Uscr_{C_d}$ of birational morphisms.

In the same way as lemma \ref{Lemme : fiber_computation}, we get
\begin{Prop}
The exceptional divisor of $\fcal_i$ is $\Dscr_i \coloneqq \Dscr \cap \Uscr_i \simeq [\C^{d-1} \times \{0\}/\Z^{n-d}]$ where the $\{0\}$ is in the $i$th position.
\end{Prop}

In the classical case, the exceptional divisor is a $\P^{n-|I|-1}$-fiber bundle over the divisor we have blown-up. In the irrational case, we study the gluing of the intersections $\mathcal{d}_i$ of the $\Dscr_i$ with $[\{0\}^I\times \C^{I^c} /\Z^{n-d}]$ i.e. the "fiber" of $\fcal_i$ over a point of $[\C^I \times \{0\}^{I^c}/\Z^{n-d}]$ (cf. Lemma \ref{Lemma : underlyingfibrespace}). We can note that the stack $\mathcal{d}_i$ is isomorphic to $[0 \times \C^{I^c \setminus \{i\}}/\Z^{n-d}]$.

It remains to us to examine the gluing of these fibers.
\\
Let $i,j$ two distinct elements of $I^c$. Then, the linear morphism $A_j^{-1}A_i$ is defined by:
\[
A_j^{-1}A_i : e_l \mapsto \begin{cases}
e_l &\text{ if } l \neq i,j \\
e_j &\text{ if } l=i \\
e_j-\sum_{m \in I^c \setminus \{j\}} \frac{\alpha_m}{\alpha_j} e_m &\text{ if } l=j
\end{cases}
\]
Note that this map fits into the following commutative diagram:

\[\begin{tikzcd}[ampersand replacement=\&]
	{\C^d} \& {\C^d} \\
	{\C^{I^c\setminus\{i\}}} \& {\C^{I^c\setminus\{j\}}}
	\arrow["{pr_{I^c\setminus\{i\}}}"', from=1-1, to=2-1]
	\arrow["{A_j^{-1}A_i}", from=1-1, to=1-2]
	\arrow["{\widetilde{A_{ji}}}"', from=2-1, to=2-2]
	\arrow["{pr_{I^c\setminus\{j\}}}", from=1-2, to=2-2]
\end{tikzcd}\]
The morphism $\widetilde{A_{ji}}$ is an isomorphism defined by:  
\[
\widetilde{A_{ji}} : e_l \mapsto \begin{cases}
e_l &\text{ if } l \in I^c \setminus\{ i,j\} \\
-\sum_{m \in I^c \setminus \{j\}} \frac{\alpha_m}{\alpha_j} e_m &\text{ if } l=j
\end{cases}
\]
and its inverse is 
\[
\widetilde{A_{ji}}^{-1} : e_l \mapsto \begin{cases}
e_l &\text{ if } l \in I^c \setminus\{ i,j\} \\
-\sum_{m \in I^c \setminus \{i\}} \frac{\alpha_m}{\alpha_i} e_j &\text{ if } l=i
\end{cases}
\]
Note that the vector $\widetilde{A_{ji}}^{-1}(e_i)$  do not depend on $j$. We deduce the following statement:

\begin{Prop}
The gluing of the $\mathcal{d}_i$ is the quantum projective space of dimension $n-|I|-1$ described by the vectors
\[e_j, j \in I^c \setminus \{i\}, v_d=-\sum_{m \in I^c \setminus i } \frac{\alpha_m}{\alpha_i} e_m\] and the calibration $pr_{I^c \setminus \{i\}}A_i^{-1}hH_i : \Z^n \to pr_{I^c \setminus \{i\}}(A_i^{-1} \Gamma)$ (with $\Ical \setminus \{k\}$ as virtual generators set) where $(A_i,H_i)$ is the pair of linear morphisms associated to the cone $\sigma_i$ for a fixed $i \in I^c$.
\end{Prop}

\begin{Rem} \label{Rem : nonsimp_blowup}
    This construction works for any simplicial cone. If we consider instead a non-simplicial cone $\sigma$, we have to take $\alpha$ in the interior of one face $\tau$ of $\sigma$ and the blow-up of $\Uscr_\sigma$ is the quantum toric stack given by the star subdivision $\sigma^*(\alpha)$ of $\sigma$ (see \cite[section 11.1]{cox}) i.e. the fan containing the following cones
   \begin{itemize}
    \item $\theta \preceq \sigma$ if $\alpha \notin \theta$ ; 
    \item $\Cone(\theta,\alpha)$ if  $\alpha \notin \theta \preceq \sigma$ and ${\alpha} \cup \theta \subset \tau \preceq  \sigma$
\end{itemize}
\end{Rem}

\section{Decomposition of birational toric morphisms by blow-ups and blow-downs}
\label{9-Decomp_Def : morph_birat}

A classic result in birational geometry is that a birational morphism between two proper smooth algebraic varieties (over an algebraically closed field) decomposes into a sequence of blow-ups and blow-downs of smooth centers (cf. \cite{wlodarczyk1999birational} and \cite{abramovich2002torification}). The proof of this theorem comes from \nom{Wlodarczyk}'s proof of the toric case. The statement we are interested in here is the following slightly weaker result:
 
    \begin{Thm}[\cite{10.2307/2155318} Theorem A] \label{Thm : decomp_rat_blowup}
Let $\Delta,\Delta'$ be two simplicial fans of $\R^n$ such that $|\Delta|=|\Delta'|$. Then there is a sequence of simplicial fans $(\Delta_i)_{i=0..n} $ such that $\Delta=\Delta_0, \Delta_n= \Delta'$ and each $\Delta_i$ is obtained from $\Delta_{i-1}$ by a blow-up or a blow-down, $i=1,\ldots,n$.
\end{Thm}

The analogous result for the quantum case is as follows: 

\begin{Thm} \label{decomp_blowup}
Let $\Delta, \Delta'$ two families of strongly convex polyhedral simplicial cones of $\R^d$ stable by intersection and by taking faces. If $|\Delta|=|\Delta'|$ then there exists two calibrations $(h,\Ical)$, $(h,\Ical')$\footnote{This is the same morphism $h$ in the two calibrations} such that $(\Delta,h,\Ical)$ and $(\Delta',h,\Ical')$ are quantum fans and such that there exists a sequence of simplicial quantum fans $(\Delta_i)_{0 \leq i \leq n}$  such that $\Delta=\Delta_0, \Delta_n= \Delta'$ and $\Delta_i$ is obtained from $\Delta_{i-1}$ by blow-up or blow-down for every $i=1,\ldots,n$.
\end{Thm}

\begin{proof}

Since the cones of the two families are simplicial, then by \cite[Lemma 2.2.4]{boivin2023moduli}, there exist $\Delta_{rat}$ and $\Delta'_{rat}$ composed of rational cones and "close" $\Delta$ and $\Delta'$ (for any metric on $\R^d$  on generators of the cones). We apply theorem \ref{Thm : decomp_rat_blowup} to connect $\Delta_{rat}$ and $\Delta'_{rat}$ by blow-ups and blow-downs. Let $G$ be the set of cone generators of all $\Delta_i$ and let $h_{rat} : \Z^G \to \R^d$ the epimorphism defined by these generators. Let $h$ be an epimorphism $\Z^G \to \R^d$ "close" to $h_{rat}$ and such that there exist $I$ and $I'$ subsets of $G$ such that for any $i \in I$ (resp. $i \in I'$), $h(e_i)$ is a generator of a 1-cone of $\Delta$ (resp. $\Delta'$) and all generators are of this form (this also exists thanks to the openness of the spaces of $D$-admissible calibration $\Omega(D)$ , again by \cite[Lemma 2.2.4]{boivin2023moduli}). Finally, we pose $\Ical'\coloneqq G \setminus I$ and $\Ical'\coloneqq G \setminus I'$ to complete this proof.

\end{proof}

\begin{Ex} 
Let $(\Delta,h : \Z^n \to \Gamma \subset \R^d,\Ical)$ be a complete quantum fan such that there exists $k \in \Ical$ such that $h(e_k) \in \R_{<0}^d$. \\
Then the toric stack $\Xscr_{\Delta,h,\Ical}$ can be linked, by blow-ups and blow-downs, to the quantum projective space given by the fan $(\Delta',h,\Ical \setminus \{k\} \cup (\Delta(1) \setminus \{1,\ldots,d\}))$ where $\Delta'$ is the fan whose maximum cones are 
\[
\Cone(e_1,\ldots,e_d),\Cone(e_1,\ldots,\widehat{e_i},h(e_k),\ldots,e_d),\ldots, \Cone(e_1,\ldots,e_{d-1},h(e_k))
\]
More precisely, we first blow-up the substack corresponding to a cone containing $h(e_k)$ and then blow-down all the other $1$-cones except $e_1,\ldots,e_d$ to obtain the quantum projective space (see figure \ref{fig:bir_morph_from_proj_sp}).

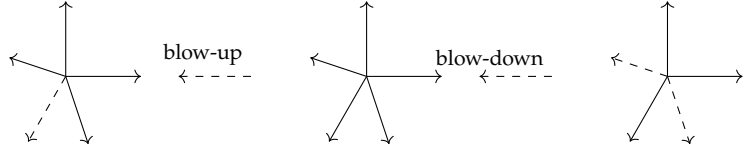
\begin{figure}[!ht]
    \centering
     \begin{tikzpicture}[scale=0.5]
\draw[->] (0,0) -- (0,2); 
\draw[->] (0,0) -- (2,0);
\draw[->] (0,0) -- (-1.5,.5);
\draw[dashed,->] (0,0) -- (-1,-1.738);
\draw[->] (0,0) -- (0.6,-1.8) ;

\draw[<-, dashed] (3,0)--(5,0) node[above left]{\tiny blow-up};
\draw[<-, dashed] (11,0)--(13,0) node[above left]{\tiny blow-down};

\draw[->] (8,0) -- (8,2); 
\draw[->] (8,0) -- (10,0);
\draw[->] (8,0) -- (6.5,.5);
\draw[->] (8,0) -- (7,-1.738);
\draw[->] (8,0) -- (8.6,-1.8) ;

\draw[->] (16,0) -- (16,2); 
\draw[->] (16,0) -- (18,0);
\draw[->,dashed] (16,0) -- (14.5,.5);
\draw[->] (16,0) -- (15,-1.738);
\draw[->, dashed] (16,0) -- (16.6,-1.8) ;

\end{tikzpicture}
    \caption{Birational morphism from a quantum projective space}
    \label{fig:bir_morph_from_proj_sp}
\end{figure}

\end{Ex}

\begin{Cor} \label{raffinement}
Let $(\Delta,h,\Ical)$ be a quantum fan and $\Delta'$ be a refinement of $\Delta$ such that $\Delta'(1) \setminus \Delta(1) \subset \Ical$. Then the natural morphism $\Xscr_{\Delta,h,\Ical} \dashrightarrow \Xscr_{\Delta',h,\Ical}$ is a composition of blow-ups and blow-downs.
\end{Cor}

\section{Cobordisms}
\label{chap:10} 
This section is devoted to the study of birational morphisms described by combinatorial constructions that are more elaborate than star subdivisions. More precisely, we will define a generalization of polytope cobordisms, described in particular in \cite{bosio2006}, to the general case of quantum fans. In the same way that these cobordisms allow us to obtain birational maps between LVM manifolds and, under certain conditions of rationality, between projective toric orbifolds induced by quotients, fan cobordisms allow us to obtain birational morphisms between quantum toric stacks that have different combinatorial types.

\subsection{Rational case: LVM manifolds and surgery}
\label{10-Polytopes_quadriques} 
\subsubsection{Polytopes and real quadrics intersection }

Recall the classical framework described in \cite{bosio2006} : \\
Let $A=(A_1,\ldots,A_n) \in \Mcal_{n,p}(\R)$. The matrix $A$ is associated with the system: 
\begin{align} 
     \begin{cases}
         \displaystyle \sum_{i=1}^n A_i|z_i|^2=0 \\
   \displaystyle \sum_{i=1}^n |z_i|^2=1
     \end{cases}
\end{align}

Note $X_A$ the set of its solution in $\C^n$ (it is an intersection of real quadrics in the sphere $\S^{2n-1}$). This is (see \cite[Lemma 1.1]{LV}) a real manifold (of dimension $2n-p-1$) if, and only if, 
\begin{itemize}
    \item $0_{\R^p} \in \Conv(A_1,\ldots,A_n)$ [Siegel condition]; 
    \item If a set $I \subset \{1,\ldots,n\}$ satisfies $0 \in \Conv(A_i, i \in I)$ then $|I|>p$ [Weak hyperbolicity].
\end{itemize}

In what follows, we suppose that $A$ satisfies these conditions (we say that $A$ is admissible).

The torus $(\S^1)^n$ has a naturel action on $X_A$ (through the inclusion $(\S^1)^n \hookrightarrow (\C^*)^n$).

\begin{Prop}[\cite{bosio2006} lemma 0.12]
The quotient $X_A/(\S^1)^{n}$ is a simple polytope of dimension $n-p-1$. 
\end{Prop}
We note $P_A$ this polytope.

\begin{Ex} \label{ex_LVM1} 
Take $A=\begin{pmatrix} 1 & 0 & -1 & -1 &-1 \\ 0 & 1 & -1 & -1 & -1 \end{pmatrix}$. \\
Then $X_A$ is the set of solutions of 
\[
\begin{cases}
|z_1|^2=|z_3|^2+|z_4|^2+|z_5|^2 \\
|z_2|^2=|z_3|^2+|z_4|^2+|z_5|^2 \\
|z_1|^2+|z_2|^2+|z_3|^2+|z_4|^2+|z_5|^2=1. \\
\end{cases}
\]
Hence, $X_A=\frac{1}{\sqrt{3}}(\S^1 \times \S^1 \times \S^5)$. The polytope $P_A$ is 
\[
\S^1 \times \S^1 \times \S^5/ \S^1  \times\S^1\times (\S^1)^3=\S^5 /(\S^1)^3=\P^2/(\S^1 \times \S^1)= \text{full triangle} 
\]
\end{Ex}

\begin{Ex}[\cite{LV} Example 5.6]\label{ex_LVM2}
Take $A=\begin{pmatrix} 1 & 0 & 5 & 1 &-2 \\ 0 & 1 & 3 & 0 & -2 \end{pmatrix}$. \\
Then  $X_A \simeq \S^3 \times \S^3 \times \S^1$ and the polytope $P_A$ is a cut full triangle.
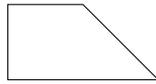
\begin{figure}[!ht]
    \centering
\begin{tikzpicture}
\draw (0,0) -- (2,0) -- (1,1) -- (0,1) -- (0,0) ; 
\end{tikzpicture}
    \caption{Cut full triangle}
    \label{fig:triangle_coupé}
\end{figure}

\end{Ex}

Every simple polytope can be realized in this way:

\begin{Thm}[\cite{bosio2006} Theorem 0.14]
Let $P$ be a simple polytope. There exists $n,p \in \N_{\geq 0}$ and a matrix $A \in \Mcal_{n,p}(\R)$ such that $P=P_A$.
\end{Thm}

Restrict the natural action of $(\S^1)^n$ on $X_A$ to an action of $\S^1$ through the diagonal morphism
\[z \in \S^1 \mapsto (z,\ldots,z) \in (\S^1)^n.\]
The quotient $\widetilde{X}_A$ of $X_A$ by this action is
a real manifold described by
\[
\widetilde{X}_A=\left\{[z_1:\ldots:z_n] \in \C\P^{n-1} \mid \sum_{i=1}^n A_i |z_i|^2=0\right\}
\]

\begin{Thm}[\cite{bosio2006} Theorem 12.2]
Let $A$ be an admissible matrix of $\Mcal_{n,p}(\R)$. Then there are two possible cases:
\begin{itemize}
    \item If $p$ is odd, then the manifold $X_A$ can be endowed with a complex structure;
    \item If $p$ is even, then the manifolds $\widetilde{X}_A$ and $X_A \times \S^1$ can be endowed with a complex structure.
\end{itemize}
\end{Thm}

\begin{Ex}
In example \ref{ex_LVM1}, the manifold $X_A/\S^1=\S^5 \times \S^1$ can be endowed with the complex structure of a Hopf manifold (see \cite{Hopf}) and in example \ref{ex_LVM2}, the manifold $X_A/\S^1=\S^3 \times \S^3$ can be endowed with the complex structure of a Calabi-Eckmann manifold (see \cite{10.2307/1969750}). . \\ The manifold $X_A \times \S^1$ is the product of the previous manifolds with a complex (compact) torus.
\end{Ex}

In both cases, the complex manifolds obtained are LVM manifolds. These manifolds are non-Kähler compact complex manifolds fully combinatorially described by a family of points $\Lambda_i$ in $\C^m$ satisfying the Siegel condition and the weak hyperbolicity. They were introduced by S. \nom{Lopez de Medrano}, A.\nom{Verjovsky} (in \cite{de1997new} for $m=1$) then generalized by L.\nom{Meersseman} (in \cite{meers1998} for any $m$). Here, we have $m=\frac{p+1}{2}$ for odd $p$ or $\frac{p}{2}$ or $\frac{p+2}{2}$ if $p$ is even (depending on the case). Under rationality conditions on $\Lambda_i$, we can prove that LVM manifolds are Seifert fiber bundles in complex tori (of complex dimension $m$) over a projective toric orbifold (see \cite[Theorem A]{LV}). More generally, it can be proved, under an assumption about virtual generators, that LVM manifolds are $\C^m$-fiber bundles over a toric quantum stack (see \cite[Theorem 9.13]{katzarkov2020quantum}). We recover the previous point of view since, under these rationality conditions, the toric orbifold is a $\Z^{2m}$-gerbe over the quantum toric stack (we thus recover the $\C^m/\Z^{2m}=(\S^1)^{2m}$-fiber bundle).

\subsubsection{Polytope cobordisms and surgery}
\label{10-cob_polytopes}

In this subsubsection, we see how operations on polytopes can describe birational transformations between the associated manifolds of the previous subsection (see \cite{bosio2006}, \cite[section 3]{timorin1999analogue}, \cite{McMullen1993}).

\begin{Def}[\cite{bosio2006} definition 2.1, \cite{timorin1999analogue}]
\label{cobord_polytope}
Let $P$ and $Q$ be two simple polytopes of the same dimension $q$. Let $W$ be a polytope of dimension $q+1$. We say that $W$ is a cobordism between $P$ and $Q$ if $P$ and $Q$ are two disjoint facets of $W$. We say that $W$ is a trivial (resp. elementary) cobordism if $W \setminus (P \cup Q)$ contains no (resp. one) vertex of $W$. 
\end{Def}

Note that a cobordism is the concatenation of elementary cobordisms. We can therefore restrict ourselves to the study of elementary cobordisms.

Consider an elementary cobordism $W$ between two polytopes $P$ and $Q$. Let $v$ be the vertex of $W \setminus (P \cup Q)$. By uniqueness, every edge starting from $v$ is connected to a vertex of $P$ or a vertex of $Q$. Let $a$ be the number of edges connecting $v$ to a vertex of $P$ and $b$ the number of edges connecting $v$ to a vertex of $Q$. By simplicity of $W$, we have $a+b=\dim(W)=q+1$. 

\begin{Def}[\cite{bosio2006} Definition 2.2, \cite{timorin1999analogue} 3.1]
The index of the cobordism $W$ is the pair of natural number $(a,b)$. We say that $Q$ is obtained by a flip of $P$ of index $(a,b)$. 
\end{Def}
The cobordism relation being transitive, if $Q$ is obtained by a flip of $P$ of index $(a,b)$ then $Q$ is obtained by a flip of $P$ of index $(b,a)$. A cobordism can be seen as a continuous deformation of the polytope $P$ and $Q$. More precisely, the simplex given by the $a$ vertices of $P$ connected to $v$ is reduced to a point and then expanded into the simplex formed by the $b$ other points connected to $v$ to obtain $Q$ (as described before the definition 2.4 of \cite{bosio2006}).

\begin{Ex}

\label{ex_cobord}
We have a cobordism between the polytopes of examples \ref{ex_LVM1} and \ref{ex_LVM2} given by the polytope in figure \ref{cob1}.\\

\begin{figure}[!ht]
    \centering
    \begin{tikzpicture}
\draw (0,0) -- (1,0) -- (0.5,.5) ;
\draw (0,0)--(0,.5) ; 
\draw[dashed] (1,1) -- (2,1) ; 
\draw (2,1) -- (1,2)  ;
\draw[dashed] (1,2) -- (1,1) ;
\draw[dashed] (0,0) -- (1,1);
\draw (1,0) -- (2,1);
\draw (0.5,1.5)-- (.5,.5) -- (0,.5) -- (.5,1.5) ; 
\draw (1,2)--(.5,1.5) ;
\end{tikzpicture}
    \caption{Cobordism between the triangle and the cut triangle}
    \label{cob1}
\end{figure}
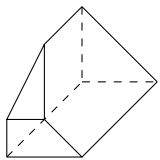

\end{Ex}

\begin{Def}
The polytope obtained after reducing the points connected to the vertex of $W$ is called the catastrophe polytope of the elementary cobordism $W$. We note $C$ this polytope.
\end{Def}
\begin{Rem}
If $a=1$ then $C=P$ and if $b=1$ then $C=Q$
\end{Rem}
\begin{Avert} \label{Warning : nonsimple}
The fact that $P$ and $Q$ are simple does not imply that $C$ is. More precisely, if $a>1$ and $b>1$ then $C$ is not simple (see example \ref{ex : cob_cube}).
\end{Avert}

Another way of looking at it is that the transition from $P$ to $Q$ is achieved by “cutting out” a neighborhood of the face generated by the $a$ points of $W$ connected to $P$ and by gluing back together a closed neighborhood of the face generated by the $b$ points of $Q$.

\begin{Def} \label{Def : polyt_transition}
The polytope obtained by cutting either $P$ or $Q$ in such way is called the transition polytope of the elementary cobordism $W$.
\end{Def}

\begin{Rem}
If $a=1$ then $T=Q$ and if $b=1$ then $T=P$
\end{Rem}

\begin{Ex} \label{ex : cob_cube}
We can consider a cobordism of index (2,2) from the cube (the final state is pictured in figure \ref{fig:cube et livre pentagonal})  : 

\begin{figure}[!ht]
    \centering
    \begin{tikzpicture}
\draw (0,0) -- (1,0) -- (1,1) -- (0,1) -- (0,0);
\draw[dashed] (0.6,0.6) -- (1.6,0.6) ; 
\draw (1.6,0.6)-- (1.6,1.6) -- (0.6,1.6)  ;
\draw[dashed] (0.6,1.6) -- (0.6,0.6) ;
\draw[dashed] (0,0) -- (0.6,0.6);
\draw (1,0) -- (1.6,0.6);
\draw (0,1) -- (0.6,1.6);
\draw (1,1) -- (1.6,1.6);

\draw (6,0) -- (6.5,0) --  (6,.5) -- (6,0);
\draw[dashed] (6.6,0.6) -- (7.6,0.6) ; 
\draw (7.6,0.6)-- (7.6,1.1) -- (6.6,1.6)  ;
\draw[dashed] (6.6,1.6) --(6.6,0.6) ;
\draw[dashed] (6,0) --(6.6,0.6);
\draw (7.3,0.3) -- (7.6,0.6);
\draw (6,.5) -- (6.6,1.6);
\draw (7.3,0.3) -- (7.6,1.1);
\draw (6.5,0) -- (7.3,.3);
\end{tikzpicture}
    \caption{Initial and final polytope}
    \label{fig:cube et livre pentagonal}
\end{figure}
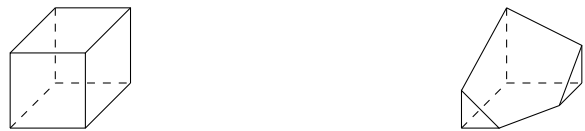

One edge of the cube has been reduced to a point and one edge of the final state has been reduced to a point, but in the “other direction”. The transition polytope is shown in figure \ref{fig:transipol}.  

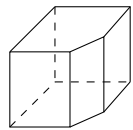
\begin{figure}[!ht]
    \centering
  \begin{tikzpicture}
\draw (0,1)--(0,0) -- (0.8,0) -- (.8,1) -- (0,1) ; 
\draw[dashed] (0.6,0.6) -- (1.6,0.6) ; 
\draw (1.6,0.6)-- (1.6,1.6) -- (0.6,1.6)  ;
\draw[dashed] (0.6,1.6) -- (0.6,0.6) ;
\draw[dashed] (0,0) -- (0.6,0.6);
\draw (1.25,0.2) -- (1.6,0.6);
\draw (0,1) -- (0.6,1.6);
\draw (1.25,1.2) -- (1.6,1.6);
\draw (.8,1)--(1.25,1.2)--(1.25,.2) -- (.8,0); 
\end{tikzpicture}
    \caption{Transition polytope}
    \label{fig:transipol}
\end{figure}
\begin{figure}[!ht]
    \centering
\begin{tikzpicture}
\draw (0,0) -- (1,0) --  (0,.5) -- (0,0);
\draw[dashed] (0.6,0.6) -- (1.6,0.6) ; 
\draw (1.6,0.6)-- (1.6,1.1) -- (0.6,1.6)  ;
\draw[dashed] (0.6,1.6) --(0.6,0.6) ;
\draw[dashed] (0,0) --(0.6,0.6);
\draw (1,0) -- (1.6,0.6);
\draw (0,.5) -- (0.6,1.6);
\draw (1,0) -- (1.6,1.1);
\end{tikzpicture}
    \caption{Catastrophe polytope}
    \label{fig:catpol}
\end{figure}
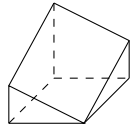

It is not simple, as it has one vertex from which four edges start.

\end{Ex}

\begin{Thm}[\cite{bosio2006} Theorem 4.8]
Let $A,B \in\Mcal_{n,p}(\R)$ two admissible matrices. Then $P_A$ and $P_B$ are linked by a flip of index $(a,b)$ (here, $a+b=n-p=\dim(P_A)+1=\dim(P_B)+1$) if, and only if, the manifolds $X_A$ and $X_B$ are linked by a surgery (see \cite[p 15]{book:91771494}) of the form 
\[
X_A \setminus \left(
\S^{2a-1} \times \D^{2b} \times (\S^1)^p
\right) \coprod_i \left(\overline{\D^{2a}}\times \S^{2b-1} \times (\S^1)^p \right)=X_B
\]
where $i $ is the $(\S^1)^n$-equivariant isomorphism defined by the composition
\[\partial \left(\S^{2a-1} \times \D^{2b} \times (\S^1)^p\right) \simeq \S^{2a-1} \times \S^{2b-1} \times (\S^1)^p  \simeq \partial \left(\overline{\D^{2n}}\times \S^{2b-1} \times (\S^1)^p\right)
\]
where $\D^n$ is the unit open ball of $\R^n$ and $\overline{\D^n}$ the unit closed ball and where the action of $(\S^1)^n=(\S^1)^a \times (\S^1)^b \times (\S^1)^p$ on $\C^a \times \C^b \times (\S^1)^p $ is given by  the inclusion $\S^1 \subset \C^*$.
If $X_A$ and $X_B$ are complex manifolds (or up to quotient by $\S^1$) then this surgery is holomorphic. 
\end{Thm}

In the classical case, these surgeries descent to quotient and induce surgeries between toric varieties (since disks are invariant by the action of compact tori of the same dimension).

This is no longer the case in the irrational case. Indeed, in this case, there is an open subset of orbits that are no longer compact (those from points where the isotropy group is not of maximal rank) and therefore the spheres' equivariant neighborhoods are no longer included in balls. \\
There is, however, a major interest in the irrational framework: we can consider the continuous family of toric quantum stacks given by the cobordism indexed by $t \in [-1,1]$ where $W_{-1}=P$ and $W_1=Q$. To study this, we will first extend these definitions to the general case of quantum fans. We use these definitions to define families in \cite[subsection 3.5]{boivin2025gluingmodulispacesquantum}.
\subsection{General case: fan cobordisms}
\label{fan_cobord}

This subsection is devoted to generalize the consideration on polytope cobordism to the general case of fans.

\begin{Prop} \label{Prop : equiv_polyt-event}
Let $W$ be a cobordism between the polytopes $P$ and $Q$ of $\R^d$. Note $i^{aff}_P : P \hookrightarrow W$ and $i^{aff}_Q : Q \hookrightarrow W$ the inclusion morphisms and $i_P$ and $i_Q$ the associated linear morphisms. Note $N_P,N_Q,N_W$ the normal fans of, respectively, $P$, $Q$ and $W$ and $N_{W,P}$, $N_{W,Q}$ the subfans of $N_W$ associated to the vertices of $P$ and $Q$ (see \cite[Chapter 2]{cox} for the details). Then,   
\begin{align}
    i_P^\top (N_{W,P})=N_P \\
    i_Q^\top(N_{W,Q})=N_Q
\end{align}

\end{Prop}

\begin{proof}
Up to affine transformation, we can suppose that the affine space generated by $i^{aff}_P(P)$ is $\{(x,t) \in \R^d \times \R\mid t=-1\}$. Then, $i_P$ is the morphism $x \in \R^d \mapsto (x,0)$ and $i_P^\perp $ is the projection $(x,t) \in \R^d \times \R \mapsto x \in \R^d$. 
We deduce that for all vertex $v$ dans $i^{aff}_P(P)$, the maximal cone $C_{v,W}$ associated to $v$ in $W$ is in form
\[
C_{v,W}=\Cone(i_P(C_{v,P}),e_{d+1})
\]
where $C_{v,P}$ is the maximal cone associated à $v$ in $P$
We get the desired equality:
\[
i_P^\top(C_{v,W})=C_{v,P}
\]
\end{proof}

Thanks to this proposition, the following definition is the transposition of the polytope cobordism definition from the previous section (cf. definition \ref{cobord_polytope}) into the fan framework (the changes come from the fact that the polytope/fan change through the normal fan flips the inclusions: the polytope vertices become maximal cones and the facets become 1-cones)

\begin{Def} \label{Def : precobordisme}
Let $(\Delta_0,h_0 : \Z^{n_0} \to \Gamma_0 \subset \R^d,\Ical_0)$, $(\Delta_1,h_1 : \Z^{n_1} \to \Gamma_1 \subset \R^d,\Ical_1)$ be two simplicial quantum fans in $\R^d$. An elementary precobordism between them is a quintuplet $(\Delta,\widetilde{\Delta_0},\widetilde{\Delta_1},\pi_0,\pi_1)$ where $\Delta=(\Delta,h : \Z^n \to \Gamma,\Ical)$ is a complete simplicial quantum fan of $\R^{d+1}$, the $\widetilde{\Delta}_i$ are subfans of $\Delta$ i.e. given by a subset of cones of $\Delta$ and the $\pi_i=(L_i : \R^{d+1} \to \R^{d},H_i : \Z^n \to \Z^{n_i})$ are  pairs of linear maps such that:
\begin{itemize}
    \item the sets $\widetilde{\Delta}_{0,max}$ and $\widetilde{\Delta}_{1,max}$ of maximal cones of, respectively, $\widetilde{\Delta}_0$ and $\widetilde{\Delta_1}$ are disjoint;
    \item $|\Delta_{max} \setminus (\widetilde{\Delta}_{0,max} \cup \widetilde{\Delta}_{1,max})|=1$
    \item The morphisms $L_i,H_i$ are onto, verify the equality $L_i h_i=h H_i$, $L_i$ sends the cone of $\widetilde{\Delta_i}$ on the cones of $\Delta_i$ and $H_i$ sends the cones of $\widetilde{\Delta_{i,h}}$ on the cones of $\Delta_{i,h_i}$
\end{itemize}
\end{Def}

This is only a \textit{pre}-cobordism, as it does not take into account the full combinatorial data of a quantum fan. Indeed, the number of generators $n_0,n_1$ and $n$, and the different sets of virtual generators $\Ical_0,\Ical_1$ and $\Ical_2$ have no relation to each other in this definition. So we are have to rigidify this definition. To do this, we give some properties of pre-cobordisms in order to add some conditions to definition \ref{Def : precobordisme}.

Let $(\Delta,\widetilde{\Delta_0},\widetilde{\Delta_1},\pi_0,\pi_1)$ be a pre-cobordism between the quantum fans $(\Delta_0,h_0 : \Z^{n_0} \to \Gamma_0,\Ical_0)$ and $(\Delta_1,h_1 : \Z^{n_1} \to \Gamma_1,\Ical_1)$. \\
Note $\sigma$ the unique maximal cone in $\Delta_{max} \setminus (\widetilde{\Delta}_{0,max} \cup \widetilde{\Delta}_{1,max})$. Note $\sigma_1,\ldots,\sigma_{d+1}$ the maximal cones of $\Delta$ which intersect $\sigma$ in a face of dimension $d$. 
 
Then, note $a$ the number of $\sigma_i$ such that $\sigma_i \in \widetilde{\Delta_0}$ and $b=d+1-a$ the number of $\sigma_i$ in $\widetilde{\Delta_1}$. 

\begin{Def}
The couple $(a,b)$ is the index of the pre-cobordism $(\Delta,\widetilde{\Delta_0},\widetilde{\Delta_1},\pi_0,\pi_1)$.
\end{Def}

\begin{Lemme} \label{Lemme : nb_gen_cobord}
Let $(\Delta,\widetilde{\Delta_0},\widetilde{\Delta_1},\pi_0,\pi_1)$ be a pre-cobordism of index $(a,b)$ between the quantum fan $(\Delta_0,h_0 : \Z^{n_0} \to \Gamma_0,\Ical_0)$ and $(\Delta_1,h_1 : \Z^{n_1} \to \Gamma_1,\Ical_1)$. Then
\begin{itemize}
    \item if $a,b \geq 2$ then $\Delta_0(1)=\Delta_1(1)$ and 
    \[
    |\Delta(1)|=|\Delta_0(1)|+2=|\Delta_1(1)|+2,
    \]
    \item If $a=1$ then $\Delta_0(1)\varsubsetneq\Delta_1(1)$ and
    \[
    |\Delta(1)|=|\Delta_0(1)|+3=|\Delta_1(1)|+2,
    \]
    \item If $b=1$ then $\Delta_0(1)\varsupsetneq\Delta_1(1)$ and 
    \[
    |\Delta(1)|=|\Delta_0(1)|+2=|\Delta_1(1)|+3.
    \]
\end{itemize}
\end{Lemme}

\begin{proof}
Proposition \ref{Prop : equiv_polyt-event} and the duality between the combinatorics of a polytope and that of its normal fan reduce this statement to \cite[Theorem 3.4.1]{timorin1999analogue}
and \cite[Proposition 2.13]{bosio2006}.
\end{proof}

\begin{Def} \label{cobordisme_def}
Let $(\Delta_0,h_0 : \Z^{n_0} \to \Gamma_0 \subset \R^d,\Ical_0)$, $(\Delta_1,h_1 : \Z^{n_1} \to \Gamma_1 \subset \R^d,\Ical_1)$ be two simplicial quantum fans in $\R^d$. A cobordism is a pre-cobordism $((\Delta,h,\Ical),\widetilde{\Delta_0},\widetilde{\Delta_1},\pi_0,\pi_1)$ (we note $(a,b)$ its index) such that:
\begin{itemize}
    \item $n_0+2=n_1+2=n$ ; 
    \item the morphisms $H_i$ are the projection $\Z^{n} \to \Z^{n-2}$ removing the two coordinates indexed by $\Delta(1) \setminus (\Delta_0(1) \cup \Delta_1(1))$ (which is of cardinal 2 thanks to lemma \ref{Lemme : nb_gen_cobord}); 
    \item if $a=1$ (resp. $b=1$) then we note $\widetilde{\Ical_0}=\Ical \cup (\Delta_1(1) \setminus \Delta_0(1))$ (resp. $\widetilde{\Ical_1}=\Ical \cup (\Delta_0(1) \setminus \Delta_1(1) )$)  and otherwise,  $\widetilde{\Ical_0}=\Ical$ (resp. $\widetilde{\Ical_1}=\Ical$). The pairs of morphisms $(L_i,H_i)$ are quantum fan morphisms $(\widetilde{\Delta}_i,h,\widetilde{\Ical_i}) \to (\Delta_i,h_i,\Ical_i)$.
\end{itemize}
\end{Def}

In the following, the first three elements of the quintuplet $(\Delta,\widetilde{\Delta_0},\widetilde{\Delta_1},\pi_0, \pi_1)$ denote the quantum fans considered in definition \ref{Def : precobordisme} and the $\pi_i$ denote the quantum fan morphisms $(\widetilde{\Delta_i},h,\widetilde{\Ical_i}) \to (\Delta_i,h_i,\Ical_i)$.

\begin{Avert}
The couple $(id : \Z^{n+2} \to \Z^{n+2},id : \Z^{d+1} \to \Z^{d+1})$ is not, in general, a fan morphism between $(\widetilde{\Delta_i},h,\widetilde{\Ical}_i)$ and $(\Delta,h,\Ical)$ due to the change of sets of virtual generators.
 \end{Avert}
 This is the only obstruction to it being a morphism: 
 \begin{Lemme}
The couple $(id : \Z^{n+2} \to \Z^{n+2},id : \Z^{d+1} \to \Z^{d+1})$  is a toric morphism $(\widetilde{\Delta_0},h,\widetilde{\Ical}_0) \to (\Delta,h,\Ical)$ (resp. $(\widetilde{\Delta_1},h,\widetilde{\Ical}_1) \to (\Delta,h,\Ical)$) if, and only if, $a \geq 2$ (resp. $b \geq 2$).
 \end{Lemme}
 \begin{proof}
If $a=1$ then it's not a morphism because it sends a non-virtual generator to a virtual generator. If $a \geq 2$ then since $\widetilde{\Delta_0}$ (resp. $\widetilde{\Delta_0}$) is a sub-fan of $\Delta$ and $\Ical_0=\Ical$ (resp. $\Ical_1=\Ical$) then it is indeed a morphism.
 \end{proof}
 
 In the case $a=1$, the morphism we obtain between $(\Delta_0,h,\Ical_0)$ and $(\Delta,h,\Ical)$ is a birational morphism given by the composition of the inclusion morphism $(\Delta_0,h,\Ical)\hookrightarrow (\Delta, h,\Ical)$ and the birational morphism $(\Delta_0,h,\Ical_0) \dashrightarrow (\Delta_0,h,\Ical)$ rendering virtual the single element of $\Delta_1(1) \setminus \Delta_0(1)$ (not appearing in $\Ical$).

\begin{Ex} \label{ex_BlP2}

The following example is the fan translation of example \ref{ex_cobord}.\\
Let $\Delta$ be the fan of $\R^3$ whose maximal cones are :
\begin{align*}
    &\Cone(e_1, e_2,e_3), \Cone(-e_1-e_2, e_2,e_3),\Cone(e_1, -e_1-e_2,e_3), \\
    &\Cone(e_1,-e_1-e_2,-e_2-e_3),\Cone(e_1,-e_3,-e_2-e_3) ,\Cone(e_1, e_2,-e_3), \\ 
    &\Cone(-e_1-e_2, e_2,-e_3),\Cone(-e_1-e_2, -e_2-e_3,-e_3), 
\end{align*}
i.e. the fan of the blow-up of a point in $\P^2 \times \P^1$ since
\begin{equation} \label{equalite_blowup}
    -e_2-e_3=(-e_1-e_2)+e_1+(-e_3)
\end{equation}

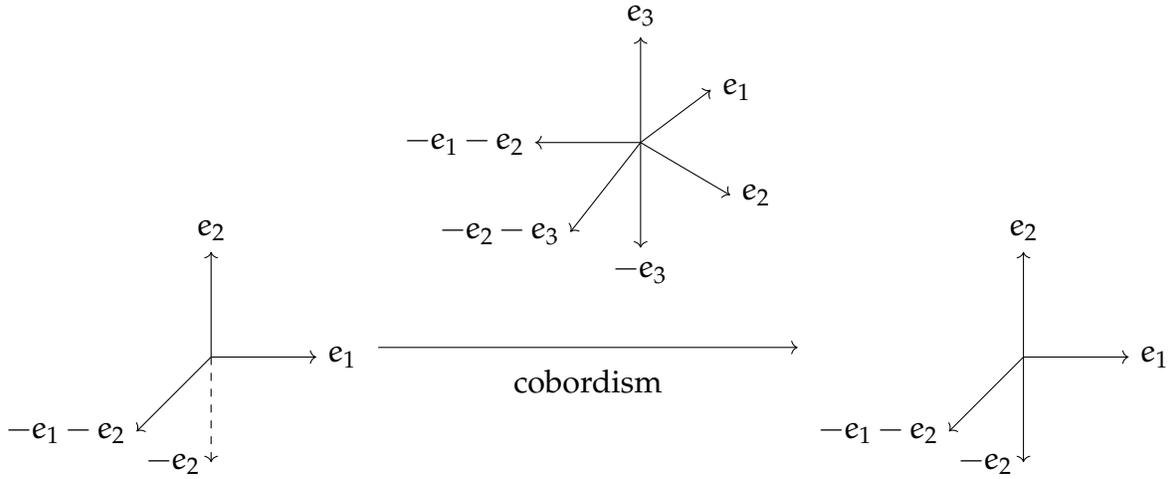
\begin{figure}[!ht]
    \centering
   
\begin{tikzpicture}[scale=.9]
\node (cob) at (0,0)
{
\begin{tikzpicture}[scale=0.7]
\draw[->] (0,0) -- (0,-2) node[below]{$-e_3$} ;
\draw[->] (0,0) -- (0,2) node[above]{$e_3$} ; 
\draw[->] (0,0) -- (-2,0) node[left]{$-e_1-e_2$} ;
\draw[->] (0,0) -- (-4/3,-1.7) node[below,left]{$-e_2-e_3$} ;
\draw[->] (0,0) -- (4/3,1) node[above,right]{$e_1$} ;
\draw[->] (0,0) -- (1.7,-1) node[below,right]{$e_2$} ;
\end{tikzpicture}
};

\node (P2) at (-6,-3)
{
\begin{tikzpicture}[scale=0.7]

\draw[->] (0,0) -- (0,2) node[above]{$e_2$} ; 
\draw[->] (0,0) -- (2,0) node[right]{$e_1$} ;
\draw[->] (0,0) -- (-1.414,-1.414) node[left]{$-e_1-e_2$} ;
\draw[dashed,->] (0,0) -- (0,-2) node[left]{$-e_2$} ;
\end{tikzpicture}
} ;
\node (P1P1) at (6,-3)
{
\begin{tikzpicture}[scale=0.7]

\draw[->] (0,0) -- (0,2) node[above]{$e_2$} ; 
\draw[->] (0,0) -- (2,0) node[right]{$e_1$} ;
\draw[->] (0,0) -- (-1.414,-1.414) node[left]{$-e_1-e_2$} ;
\draw[->] (0,0) -- (0,-2) node[left]{$-e_2$} ;
\end{tikzpicture}
};

\draw[->] (P2.east) -- (P1P1.west) node at (0,-3.5) {cobordism};

\end{tikzpicture}

 \caption{Cobordism between the fan of $\P^2$ and its blow-up in one point}
    \label{fig:cobordismeP2}
\end{figure}

Let $(h : \Z^6 \to \Z^3,\emptyset)$ be the standard calibration such that $h(e_4)=-e_2-e_3$, $h(e_5)=e_3$, $h(e_6)=-e_3$. \\
Note $\widetilde{\Delta_0}$ the subfan given by the three first cones of $\Delta$ and $\widetilde{\Delta_1}$ the four last one.\\
Let $\Delta_0$ be the fan of $\R^2$ whose maximal cones are:
\[
\Cone(e_1,e_2), \Cone(e_2,-e_1-e_2),  \Cone(-e_1-e_2,e_1)
\]
i.e. the fan of $\P^2$. \\
Let $(h_0 : \Z^4 \to \Z^2, (x,y,z,t) \mapsto (x-z,y-z-t),\Ical_0=\{4\})$ be a calibration of $\Z^2$.   \\
Let $\Delta_1$ be the fan of $\R^2$ whose maximal cones are:
\[
\Cone(e_1,e_2), \Cone(e_2,-e_1-e_2), \Cone(-e_2,-e_1-e_2), \Cone(-e_1-e_2,e_1)
\]
i.e. the fan of a blow-up of $\P^2$. \\
Let $(h_1 : \Z^4 \to \Z^2, (x,y,z,t) \mapsto (x-z,y-z-t),\Ical_1=\emptyset)$ be a calibration of $\Z^2$.  
Then the fan $(\Delta,h,\Ical)$, the two fans $(\widetilde{\Delta}_0,h,\{4\})$, $(\widetilde{\Delta_1},h,\emptyset)$ and the two projections $\pi_0=(L_0 : (x,y,z) \in \R^3 \mapsto (x,y) \in \R^2, H_0 : (x,y) \in \Z^4 \times \Z^2 \mapsto x \in \Z^4)\colon (\widetilde{\Delta}_0,h,\{4\}) \to (\Delta_0,h_0,\{4\})$ and $\pi_1=(L_1 : (x,y,z) \in \R^3 \mapsto (x,y) \in \R^2, H : (x,y) \in \Z^4 \times \Z^2 \mapsto x \in \Z^4) \colon (\widetilde{\Delta}_1,h,\emptyset) \to (\Delta_0,h_0,\emptyset)$ form a cobordism between $(\Delta_0,h_0,\{4\})$ and $(\Delta_1,h_1,\emptyset)$.

\end{Ex}

\begin{Rem}
In the previous example, we could have put (potentially irrational) coefficients into the equality \eqref{equalite_blowup} to obtain a weighted blow-up of $\P^2$ (or more precisely, a quantum toric stack such that $\P^2$ is a $\Z$-gerbe over it). 
\end{Rem} 

This construction can be done in arbitrary dimension:
\begin{Const} \label{blow-up}
Let $(\Delta_0,h_0,\Ical_0)$ be a simplicial complete quantum fan in $\R^d$. 

Let $\alpha \in (\R_{> 0})^d$ and $\sigma=\Cone(v_1=h_0(e_{i_1}),\ldots,v_d=h_0(e_{i_d})) \in \Delta_{max}$ (suppose that $h_0(e_n)=\sum_{j=1}^d \alpha_i h(e_{i_j}))$.
Let $(\Delta,h,\Ical)$ be the quantum fan defined by: 
\begin{itemize}
    \item $\Delta$ is the star subdivision (of weight $\alpha$) (see remark \ref{Rem : nonsimp_blowup} for the definition) of $\Delta_0 \times (\R_{\geq 0},0,\R_{\leq 0})$ along the cone $\sigma \times \R_{\leq 0} e_{d+1}$;
    \item $h : \Z^{n+2} \to \Gamma\coloneqq \sum_{i=1}^n h(e_i)\Z+e_{d+1}\Z, h(e_i)=h_0(e_i), 1 \leq i \leq n-1,  h(e_{n})=\sum_j \alpha_i v_{i_j}-e_{d+1}$, $h(e_{n+1})=e_{d+1}$, $h(e_{n+2})=-e_{d+1}$; 
    \item $\Ical=\Ical_0 \setminus \{n\}$;
\item $\widetilde{\Delta_1}$ is the subfan of $\Delta$ whose maximal cones are:
\begin{itemize}
    \item the elements of $(\Delta_{max} \setminus \{\sigma\}) \times \R_{\leq 0}$ ; 
    \item the cones of the form $\Cone(v_i, i \in I,h(e_{n}), -e_{d+1})$, where $I$ is a subset of $[\![1,d]\!]$ of cardinal $d-2$.
\end{itemize}
\end{itemize}

Let $\Delta_1$ be the image of $\widetilde{\Delta_1}$ by the projection  $\pi : (x,t) \in \R^d \times \R \mapsto x \in \R^d$. We get the star subdivision (of weight $\alpha$) along the cone $\sigma$ i.e. it is the fan which describes the blow-up of weight $\alpha$ of $\Xscr_{\Delta_0,h_0,\Ical_0}$. \\
The quintuplet \[(\Delta,\Delta_0 \times \R_{\geq 0},\widetilde{\Delta_1}, (\pi\colon(x,t) \mapsto x,\mathrm{pr} : \Z^{n}\oplus \Z^2 \to \Z^n), (\pi : (x,t) \mapsto x,\mathrm{pr} : \Z^{n}\oplus \Z^2 \to \Z^n))\] is a cobordism between $(\Delta_0,h_0,\Ical_0)$ and $(\Delta_1,h_1\coloneqq h_0,\Ical_1 \coloneqq \Ical_0 \setminus \{n\})$, called cobordism of the $\alpha$-blow-up of $(\Delta_0,h_0,\Ical_0)$ in $\sigma$.

\end{Const}

\begin{Rem}
A blow-up of a fan is of index $(1,d)$.
\end{Rem}

The cobordism $(\Delta,\widetilde{\Delta_0},\widetilde{\Delta_1},\pi_0,\pi_1)$ between $(\Delta_0,h_0,\Ical_0)$ and $(\Delta_1,h_1, \Ical_1)$ is, through theorem \ref{Thm : equiv_cat} and definition \ref{Def : morph_birat},  equivalent to the following quantum toric stack diagram: 
\begin{equation} \label{diagr_cobord}
    \begin{tikzcd}
	{\Xscr_{\widetilde{\Delta_0},h,\widetilde{\Ical}_0}} & {} & {\Xscr_{\Delta,h,\Ical}} & {} & {\Xscr_{\widetilde{\Delta_0},h,\widetilde{\Ical}_1}} \\
	{\Xscr_{\Delta_0,h_0,\Ical_0}} &&&& {\Xscr_{\Delta_1,h_1,\Ical_1}}
	\arrow[dashed,hook, from=1-1, to=1-3]
	\arrow[dashed,hook', from=1-5, to=1-3]
	\arrow[two heads, from=1-1, to=2-1]
	\arrow[two heads, from=1-5, to=2-5]
	\end{tikzcd}
\end{equation}
In this way, a cobordism links two quantum toric stacks with different combinatorial types (cf. remark \ref{Rem : modulispace}). An example is described in example \ref{ex_BlP2} between the combinatorial types $D_0=\{\{1\},\{2\},\{3\},\{1,2\},\{2,3\},\{3,4\}\}$ and $D_1=\{\{1\},\{2\},\{3\},\{4\},\{1,2\},\{2,3\},\{3,4\},\{4,1\}\}$, ordered by inclusion.

\begin{Not}
In the following, for simplicity's sake, we assume that the toric morphisms $\pi_i$ are given by the same pair of linear morphisms $(L,H)$. Then, up to linear isomorphism, we can assume that $L$ is the projection $(x,t) \in \R^{d+1} \mapsto x\in \R^d$.
\end{Not}

\begin{Rem} \label{Rem : projection}
The projection $L$ defines a continuous family of quantum fans and hence a continuous family of quantum toric stacks (which corresponds by proposition \ref{Prop : equiv_polyt-event} with taking the normal fan of the slice of the polytope defining the cobordism).
\end{Rem}

\begin{Prop} \label{cobordisme_transition}
Let $(\Delta,\widetilde{\Delta_0},\widetilde{\Delta_1},\pi_0,\pi_1)$ be a cobordism between the quantum fans $(\Delta_0,h_0,\Ical_0)$ and $(\Delta_1,h_1,\Ical_1)$. Note $\sigma$ be the maximal cone of $\Delta$ which is neither in $\Delta_0$ nor in $\Delta_1$. Let $\Sigma_i$ be the set of maximal cones of $\Delta_i$ with an intersection with $\sigma$ of maximal dimension.  
Then
\begin{equation}
    L\left(\bigcup_{\sigma \in \Sigma_0} \sigma \right)=L\left(\bigcup_{\tau \in \Sigma_1} \tau \right)
\end{equation}
and the quantum fans obtained by replacing the elements of $L(\Sigma_0)$ by the cone  $L\left(\bigcup_{\sigma \in \Sigma_0} \sigma\right)$ and the elements of $L(\Sigma_1)$ by the cone $L\left(\bigcup_{\tau \in \Sigma_1} \tau\right)$ coincide.
\end{Prop}

\begin{proof}
Let $(a,b)$ be the index of the cobordism $(\Delta,\widetilde{\Delta_0},\widetilde{\Delta_1},\pi_0,\pi_1)$. Note $v_1,\ldots,v_{d+1}$ the generators of $\sigma$. Then the cones of $\Delta_0$ intersecting $\sigma$ in a maximal ways are of the form
\[
\sigma_k\coloneqq \Cone(v_i, I \setminus \{k\},\alpha e_{d+1})
\]
where $k \in [\![1,d+1]\!]$. Similarly, the cones of $\Delta_1$ intersecting $\sigma$ with maximal dimension are of the form
\[
\tau_l= \Cone(v_i, I \setminus \{l\},\beta e_{d+1})
\]
where $l \in [\![1,d+1]\!]$ and $\alpha\beta<0$. 
If $a \geq 2$ and $b \geq 2$ then the union of the generators of $L\sigma_k$ coincides with the union of the generators of $L\tau_l$. We now prove that 
\begin{equation*}
    L\left(\bigcup_{\sigma \in \Sigma_0} \sigma \right)=L\left(\bigcup_{\tau \in \Sigma_1} \tau \right).
\end{equation*}
Let $x \in L\left(\bigcup_{\sigma \in \Sigma_0} \sigma \right)$ and $k \in [\![1,d+1]\!]$ such that $x \in L\sigma_k$. Let $\tau$ and $\tau'$ be cones of $\Sigma_1$. Up to re-indexing, we can assume they are of the form
\[L\tau=\Cone(Lv_1,\ldots,Lv_{d-1},Lv_d)\]
and
\[L\tau'=\Cone(Lv_1,\ldots,Lv_{d-1},Lv_{d+1})\]
By construction, we have 
\[
Lv_{d+1}=\sum_{i=1}^d \lambda_i Lv_i
\]
where, for all $i<d$, $\lambda_i \geq 0$ and $\lambda_d< 0$. We also have
\[
Lv_d=\sum_{i=1}^d \frac{\lambda_i}{-\lambda_d} Lv_i+\frac{1}{\lambda_d} Lv_{d+1}.
\]
Since $x \in L\sigma_k$ then 
\[
x = \sum_{i \neq k, i \leq d+1} \mu_i Lv_i 
\]
and hence
\begin{align*}
   x&=\sum_{i \neq k, i \leq d-1} (\mu_i+\mu_{d+1}\lambda_i) v_i +(\mu_d+\mu_{d+1} \lambda_d)v_d \\
   &=\sum_{i \neq k, i \leq d-1} \left(\mu_i+\mu_{d+1}\frac{\lambda_i}{-\lambda_d}\right) v_i +\left(\mu_d+\frac{\mu_{d+1}}{\lambda_d} \right)v_d 
\end{align*}

So $x \in L\tau_1$ if $\mu_d+\mu_{d+1}\lambda_d \geq 0$ and $x \in L\tau_2$ otherwise. In other words, $x \in L\left(\bigcup_{\tau \in \Sigma_1} \tau \right)$. Inclusion in the other direction is done in the same way. 

If $a=1$ or $b=1$, this cobordism is the cobordism of a weighted blow-up. The cones $\tau_1,\ldots,\tau_d$ (resp. $\sigma_1,\ldots,\sigma_d$) connected to the maximal cone form the star subdivision of the cone $\sigma_1$ (resp. $\tau_1$). We therefore deduce that 
\[
\sigma_1=\bigcup_{i=1}^d \tau_i.
\]
\end{proof}

\begin{Def}
The obtained fan is called catastrophe fan and is denoted $\Delta_{1/2}$.
\end{Def}

\begin{Ex}
In the case of blow-up, the catastrophe fan is the fan that we want to blow-up.
\end{Ex}

\begin{Rem}
In the same manner as the polytope cobordism (warning \ref{Warning : nonsimple}) is not simple, the catastrophe fan is not simplicial in general.
\end{Rem}

\begin{Rem}
It corresponds to the unique point of change of the combinatorics of the family of the remark \ref{Rem : projection}.
\end{Rem}

Since the quantum fans $(\Delta_0,h,\Ical_0)$, $(\Delta_{1/2},h,\Ical_{1/2})$ (where $\Ical_{1/2}=\Ical_0 \cap \Ical_1$) and $(\Delta_1,h,\Ical_1)$ have a common refinement then 
\begin{Prop}
The quantum toric stacks $(\Delta_0,h,\Ical_0)$, $(\Delta_{1/2},h,\Ical_{1/2})$ and $(\Delta_1,h,\Ical_1)$ are birational.
\end{Prop}

These birational morphisms can be described explicitly: \\
Note $\sigma$ the cone $L\left(\bigcup_{\sigma \in \Sigma_0} \sigma\right)$ and $\alpha \in \sigma$. 
We consider the blow-up
 of $\Delta$ along $\sigma$ in $\alpha$ (see remark \ref{Rem : nonsimp_blowup}) which we will note $\widetilde{\Delta}$. This fan corresponds to the transition polytope of definition \ref{Def : polyt_transition}.
The birational morphism $(\widetilde{\Delta},\widetilde{h},\widetilde{\Ical}) \dashrightarrow (\Delta_{1/2},h, \Ical_{1/2})$ has as exceptional divisor $\Dscr$ which is a quantum projective bundle of dimension $a-1$ over a quantum projective space of dimension $b-1$ (see \cite[Theorem 10.3.2.9]{Boivin})
thanks to the dual fact that the facet we add in the transition polytope is the product of a simplex of dimension $a-1$ and a simplex of dimension $b-1$, see \cite[Proposition 2.6]{bosio2006}. The projections $\Dscr \to \P^{a-1}_q$ and $\Dscr \to \P^{b-1}_q$ induce the blow-down $\Xscr_{(\widetilde{\Delta},\widetilde{h},\widetilde{\Ical})} \dashrightarrow \Xscr_{(\Delta_0,h,\Ical_0)}$ and $\Xscr_{(\widetilde{\Delta},\widetilde{h},\widetilde{\Ical})} \dashrightarrow \Xscr_{(\Delta_1,h,\Ical_1)}$ (which come from the simplex that was removed in the polytopal case between the transition polytope and the cobordism polytopes). We thus obtain the arrows $\Xscr_{(\Delta_0,h,\Ical_0)} \dashrightarrow \Xscr_{(\Delta_1,h,\Ical_1)}$, $\Xscr_{(\Delta_0,h,\Ical_0)} \dashrightarrow \Xscr_{(\Delta_1,h,\Ical_1)}$, $\Xscr_{(\Delta_0,h_0,\Ical_0)} \dashrightarrow \Xscr_{(\Delta_{1/2},h,\Ical_{1/2})}$ by composition. In one diagram,

\[\begin{tikzcd}
	& {\Xscr_{\widetilde{\Delta},\widetilde{h},\widetilde{\Ical}}} \\
	{\Xscr_{\Delta_{0},h_0,\Ical_0}} && {\Xscr_{\Delta_{1},h_1,\Ical_1}} \\
	& {\Xscr_{\Delta_{1/2},h,\Ical_{1/2}}}
	\arrow["{\text{blow-up}}", dashed, from=1-2, to=3-2]
	\arrow["\text{blow-down}"'{pos=0.5}, dashed, from=1-2, to=2-1]
	\arrow["\text{blow-down}"{pos=0.5}, dashed, from=1-2, to=2-3]
	\arrow[dashed, from=2-1, to=3-2]
	\arrow[dashed, from=2-3, to=3-2]
\end{tikzcd}\]

\begin{Ex}
Consider the fan counterpart of cobordism in example \ref{ex : cob_cube}. The cone $L\left(\bigcup_{\sigma \in \Sigma_0 } \sigma \right)$ is a non-simplicial cone with four sides. 

\begin{figure}[!ht]
    \centering
\begin{tikzpicture}[scale=1]
    \draw (0,0) -- (.7,1.2) ;
    \draw (0,0) -- (-.7,1.2) ; 
    \draw[dashed] (0,0) -- (.1,1.6) ; 
    \draw (0,0) -- (0,.8)  ;
    \draw (.7,1.2) -- (0.1,1.6) -- (-.7,1.2) -- (0,.8) -- (.7,1.2);
\end{tikzpicture}

\caption{The cone $L\left(\bigcup_{\sigma \in \Sigma_0 } \sigma \right)$}
    \label{fig:cata}
\end{figure}
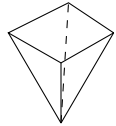
The blow-up of the cone of figure \ref{fig:cata} corresponds to its complete triangulation (see figure \ref{fig:blowupcata}).
\begin{figure}[!ht]
    \centering
\begin{tikzpicture}[scale=1]
    \draw (0,0) -- (.7,1.2) ;
    \draw (0,0) -- (-.7,1.2) ; 
    \draw[dashed] (0,0) -- (.1,1.6) ; 
    \draw (0,0) -- (0,.8)  ;
    \draw (.7,1.2) -- (0.1,1.6) -- (-.7,1.2) -- (0,.8) -- (.7,1.2);
    \draw (.7,1.2) -- (-.7,1.2) ; 
    \draw (.1,1.6) -- (0,.8)  ;
\end{tikzpicture}

\caption{The blow-up cone of $L\left(\bigcup_{\sigma \in \Sigma_0 } \sigma \right)$}
    \label{fig:blowupcata}
\end{figure}
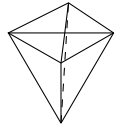
This cone can then be blow-down in two different ways:

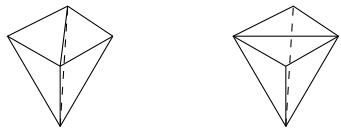
\begin{figure}[!ht]
    \centering
\begin{tikzpicture}[scale=1]
    \draw (0,0) -- (.7,1.2) ;
    \draw (0,0) -- (-.7,1.2) ; 
    \draw[dashed] (0,0) -- (.1,1.6) ; 
    \draw (0,0) -- (0,.8)  ;
    \draw (.7,1.2) -- (0.1,1.6) -- (-.7,1.2) -- (0,.8) -- (.7,1.2);
    \draw (.1,1.6) -- (0,.8)  ;
    
    \draw (3,0) -- (3.7,1.2) ;
    \draw (3,0) -- (2.3,1.2) ; 
    \draw[dashed] (3,0) -- (3.1,1.6) ; 
    \draw (3,0) -- (3,.8)  ;
    \draw (3.7,1.2) -- (3.1,1.6) -- (2.3,1.2) -- (3,.8) -- (3.7,1.2);
    \draw (3.7,1.2) -- (2.3,1.2) ; 
\end{tikzpicture}

\caption{Projections of cones of $\Delta_0$ and cones of $\Delta_1$ connected to the maximum cone of the cobordism}
    \label{fig:contracterblowupcata}
\end{figure}

The other cones remain unchanged by all these operations. The diagram of birational morphisms can therefore be written in the form  described in the figure. \ref{fig:diag_birat}. \\

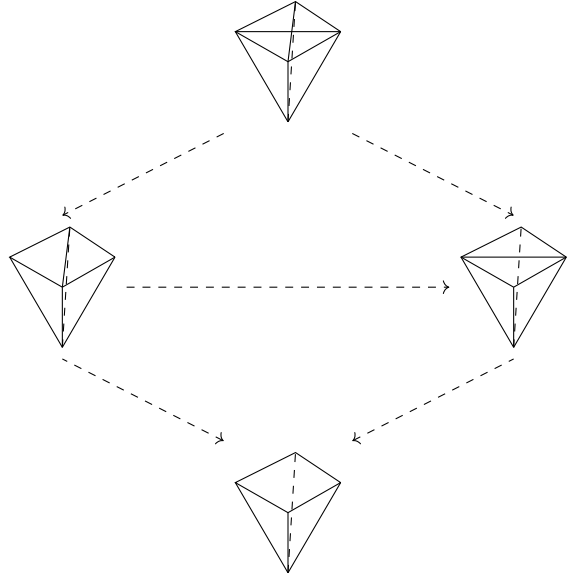
\begin{figure}[!ht]
    \centering
    
\begin{tikzpicture}

\node (bl) at (0,0)
{
\begin{tikzpicture}[scale=1]
    \draw (0,0) -- (.7,1.2) ;
    \draw (0,0) -- (-.7,1.2) ; 
    \draw[dashed] (0,0) -- (.1,1.6) ; 
    \draw (0,0) -- (0,.8)  ;
    \draw (.7,1.2) -- (0.1,1.6) -- (-.7,1.2) -- (0,.8) -- (.7,1.2);
    \draw (.7,1.2) -- (-.7,1.2) ; 
    \draw (.1,1.6) -- (0,.8)  ;
\end{tikzpicture}
};

\node (X0) at (-3,-3){

\begin{tikzpicture}
\draw (0,0) -- (.7,1.2) ;
    \draw (0,0) -- (-.7,1.2) ; 
    \draw[dashed] (0,0) -- (.1,1.6) ; 
    \draw (0,0) -- (0,.8)  ;
    \draw (.7,1.2) -- (0.1,1.6) -- (-.7,1.2) -- (0,.8) -- (.7,1.2);
    \draw (.1,1.6) -- (0,.8)  ;
\end{tikzpicture}

};

\node (X1) at (3,-3){
\begin{tikzpicture}
\draw (0,0) -- (.7,1.2) ;
    \draw (0,0) -- (-.7,1.2) ; 
    \draw[dashed] (0,0) -- (.1,1.6) ; 
    \draw (0,0) -- (0,.8)  ;
    \draw (.7,1.2) -- (0.1,1.6) -- (-.7,1.2) -- (0,.8) -- (.7,1.2);
    \draw (.7,1.2) -- (-.7,1.2) ; 
\end{tikzpicture}
};
\node (X) at (0,-6){
\begin{tikzpicture}
\draw (0,0) -- (.7,1.2) ;
    \draw (0,0) -- (-.7,1.2) ; 
    \draw[dashed] (0,0) -- (.1,1.6) ; 
    \draw (0,0) -- (0,.8)  ;
    \draw (.7,1.2) -- (0.1,1.6) -- (-.7,1.2) -- (0,.8) -- (.7,1.2);
\end{tikzpicture}
};

\draw[->,dashed] (bl.south west) -- (X0.north);
\draw[->,dashed] (bl.south east) -- (X1.north);
\draw[->,dashed] (X0.east) -- (X1.west);
\draw[<-,dashed] (X.north west) -- (X0.south);
\draw[<-,dashed] (X.north east) -- (X1.south);

\end{tikzpicture}
\caption{Diagram of birational morphisms of the cobordism}
    \label{fig:diag_birat}
\end{figure}

\end{Ex}

This diagram is reminiscent of the Atiyah flop as described in \cite[section 3]{10.2307/100556} where we blow-up the origin of the singular variety $V(xt-yz) \subset \C^4$ (which is the toric variety associated to the cone over a square). The exceptional divisor of this blow-up is $\P^1 \times \P^1$ (seen as embed into the projective space $\P^3$ by the Segre embedding). This divisor $\P^1 \times \P^1$ can be contracted onto $\P^1$ by the two projections. The two varieties thus obtained are birational and are given by the fans in figure \ref{fig:contracterblowupcata}.

To conclude this section, we look at deformations of cobordisms. In the classical framework of \cite{bosio2006}, cobordisms are easily deformed, since a (small) deformation (of vertices) of the polytope $W$ (induced by a deformation of $P$ or $Q$) does not change its combinatorics (since the realization space of simple polytopes is open).

In this way, we obtain a deformation of the associated LVM manifolds, but no deformation of the toric varieties, since we have to consider irrational coefficients.

Thanks to the quantum framework, we can easily consider deformations of cobordisms and simplicial fans since $\Omega(d,n,D)$ is an open subset of $\R^{d(n-d)}$ (cf. \cite[Lemma 2.2.4]{boivin2023moduli}: 

\begin{Thm} \label{deform_cobord} 
    Let $(\Delta,\widetilde{\Delta_0},\widetilde{\Delta_1},\pi_0,\pi_1)$ be a cobordism between the complete simplicial quantum fans $(\Delta_0,h_0,\Ical_0)$ and $(\Delta_1,h_1,\Ical_1)$. 
    Let $(\Delta'_0,h'_0,\Ical_0)$ be a quantum fan in a neighborhood of $(\Delta_0,h_0,\Ical_0)$ (where $h_0$ and $h'_0$ coincide on $\Z^{\Ical_0}$) in  $\Omega(\mathrm{comb}(\Delta_0))$. There exists a unique fan $(\Delta'_1,h'_1,\Ical_1) \in \Omega(\mathrm{comb}(\Delta_1))$ (where $h_1$ and $h'_1$ coincide on $\Z^{\Ical_1}$) and a unique cobordism $(\Delta',\widetilde{\Delta'_0},\widetilde{\Delta'_1},\pi_0,\pi_1)$ between them which minimize the distance between $h$ and $h'$. We have the same statement when we replace $0$ by $1$ and vice-versa.

    In a diagram: 

\[\begin{tikzcd}
	{(\widetilde{\Delta}'_0,h',\mathcal{I}_0')} && {(\Delta',h',\mathcal{I}')} && {(\widetilde{\Delta}'_1,h',\mathcal{I}'_1)} \\
	& {(\widetilde{\Delta}_0,h,\mathcal{I})} & {(\Delta,h,\mathcal{I})} & {(\widetilde{\Delta}_1,h,\mathcal{I})} \\
	{(\Delta'_0,h'_0,\mathcal{I}'_0)} & {(\Delta_0,h_0,\mathcal{I}_0)} && {(\Delta_1,h_1,\mathcal{I}_1)} & {(\Delta'_1,h'_1,\mathcal{I}'_1)}
	\arrow["\approx"{marking}, draw=none, from=1-1, to=2-2]
	\arrow["\approx"{marking}, draw=none, from=3-1, to=3-2]
	\arrow["\approx"{marking}, draw=none, from=1-3, to=2-3]
	\arrow["\approx"{marking}, draw=none, from=1-5, to=2-4]
	\arrow["\approx"{marking}, draw=none, from=3-5, to=3-4]
	\arrow[two heads, from=2-2, to=3-2]
	\arrow[hook, from=2-2, to=2-3]
	\arrow[hook', from=2-4, to=2-3]
	\arrow[two heads, from=2-4, to=3-4]
	\arrow[two heads, from=1-1, to=3-1]
	\arrow[hook, from=1-1, to=1-3]
	\arrow[hook', from=1-5, to=1-3]
	\arrow[two heads, from=1-5, to=3-5]
\end{tikzcd}\]

\end{Thm}

\begin{proof}
Let $(\Delta'_0,h'_0,\Ical_0)$ be a quantum fan in a neighborhood of $(\Delta_0,h_0,\Ical_0)$ and such that $h_0$ and $h'_0$ coincide on $\Z^{\Ical_0}$. The lemma \ref{Lemme : nb_gen_cobord} determines the unique candidate for $(\Delta_1,h_1,\Ical'_1)$. Prove that there is a cobordism between $(\Delta'_0,h'_0,\Ical_0)$ and $(\Delta'_1,h'_0,\Ical_1)$. \\
By continuity of projection, there exists a morphism $h' : \R^n \to \R^{d+1}$ in the neighborhood of $h$ in $\Omega(\mathrm{comb}(\Delta'))$ that minimizes the distance among the antecedents of $h_0$ (since $\Omega(\mathrm{comb}(\Delta'))$ is open and $\ker(L)$ is a convex). With the combinatorics of a cobordism fixed, we deduce the cobordism between $(\Delta'_0,h'_0,\Ical_0)$ and $(\Delta'_1,h'_0,\Ical_1)$.
\end{proof}

\begin{figure}[!ht]
    \centering
    \begin{tikzpicture}
\draw (0,0) -- (1.2,0.1) -- (0.5,.5) ;
\draw (0,0)--(0,.5) ; 
\draw[dashed] (1,1) -- (2,1) ; 
\draw (2,1) -- (1,2)  ;
\draw[dashed] (1,2) -- (1,1) ;
\draw[dashed] (0,0) -- (1,1);
\draw (1.2,0.1) -- (2,1);
\draw (0.5,1.5)-- (.5,.5) -- (0,.5) -- (.5,1.5) ; 
\draw (1,2)--(.5,1.5) ;
\end{tikzpicture}
    \caption{Deformation of the cobordism between the triangle and the cut triangle}
    \label{cob1_def}
\end{figure}
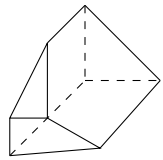 
\begin{Ex}
Let $(\Delta,\widetilde{\Delta_0},\widetilde{\Delta_1},\pi_0,\pi_1)$ be the $\alpha$-blow-up of the cone $\Cone(-e_1-e_2,e_1)$ in $(\Delta_0,h,\Ical)$ as defined in example \ref{ex_BlP2}. Then we can obtain this blow-up by a deformation of a rational blow-up (by density of $\Q^2$ in $\R^2$). 
\end{Ex}


\bibliographystyle{alpha} 
\bibliography{Biblio}

\end{document}